\documentclass{amsart}

\usepackage{amsthm}
\usepackage{amsmath}
\usepackage{amssymb}
\usepackage{amsaddr}
\usepackage{mathtools}
\usepackage{xfrac}
\usepackage{graphicx}
\usepackage{latexsym}
\usepackage{hyperref}
\usepackage{setspace}
\usepackage{enumerate}
\usepackage{bbm}
\usepackage{physics}
\usepackage{chngcntr}

 \usepackage{geometry}
 \allowdisplaybreaks

 \usepackage[dvipsnames]{xcolor}

 \usepackage{biblatex}
 \bibliography{main}
 \hypersetup{
     colorlinks=true,
     linkcolor=blue,
     filecolor=magenta,
     citecolor=green,
     urlcolor=cyan,
 }

\renewcommand{\abs}[1]{\left|#1\right|}

\newcommand{\lD}{(\lambda,D)}

\renewcommand{\epsilon}{\varepsilon}

\newcommand{\livp}[1]{$\text{(LIVP)}_{#1}$}

\newcommand{\DT}{{D}_T}
\newcommand{\Dinfty}{{D}_\infty}
\newcommand{\closedDT}{\overline{D}_T}
\newcommand{\closedDinfty}{\overline{D}_\infty}

\DeclareMathOperator{\supp}{supp}

\renewcommand{\op}[1]{\mathrm{#1}}

\newcommand{\R}{\mathbb{R}}
\newcommand{\C}{\mathbb{C}}

\newcommand{\Z}{\mathbb{Z}}
\newcommand{\N}{\mathbb{N}}
\newcommand{\on}{\quad\text{on }}
\newcommand{\xt}{{(x,t)}}

\begin{document}

\title[The evolution problem for the 1D nonlocal Fisher-KPP equation]{The evolution problem for the 1D nonlocal Fisher-KPP equation with a top hat kernel. Part 1. The Cauchy problem on the real line}

\author{D. J. Needham}
\address{School of Mathematics, University of Birmingham, 
UK, B15 2TT}

\author{J. Billingham}
\address{School of Mathematical Sciences, University of Nottingham,
    UK, NG7 2RD}

\author{N. M. Ladas}
\address{School of Mathematics, University of Birmingham, 
UK, B15 2TT}

\author{J. C. Meyer}
\address{School of Mathematics, University of Birmingham, 
UK, B15 2TT}





\begin{abstract}
We study the Cauchy problem on the real line for the nonlocal Fisher-KPP equation in one spatial dimension,
\[
u_t = D u_{xx} + u(1-\phi*u),
\]
where $\phi*u$ is a spatial convolution with the top hat kernel, $\phi(y) \equiv H\left(\frac{1}{4}-y^2\right)$. 

After observing that the problem is globally well-posed, we demonstrate that positive, spatially-periodic solutions bifurcate from the spatially-uniform steady state solution $u=1$ as the diffusivity, $D$, decreases through $\Delta_1 \approx 0.00297$ (the exact value is determined in section 3). We explicitly construct these spatially-periodic solutions as uniformly-valid asymptotic approximations for $D \ll 1$, over one wavelength, via the method of matched asymptotic expansions. These consist, at leading order, of regularly-spaced, compactly-supported regions with width of $O(1)$ where $u=O(1)$, separated by regions where $u$ is exponentially small at leading order as $D \to 0^+$. 

From numerical solutions, we find that for $D \geq \Delta_1$, permanent form travelling waves, with minimum wavespeed, $2 \sqrt{D}$, are generated, whilst for $0 < D < \Delta_1$, the wavefronts generated separate the regions where $u=0$ from a region where a steady periodic solution is created via a distinct periodic shedding mechanism acting immediately to the rear of the advancing front, with this mechanism becoming more pronounced with decreasing $D$. The structure of these transitional travelling wave forms is examined in some detail.



\end{abstract}


\maketitle

\section{Introduction}\label{section1}
Nonlocal reaction-diffusion equations arise in many different scientific areas (see, for example, \cite{Britton, Coombes,coville_dupaigne_2007, Furter1989,Surulescu,  Volpert2014} and \cite{Kavallarisbook,Volpert2009}). Many of these applications are biomedical, including tumour modelling and models for evolution and speciation (see \cite{VolpBio} for an extensive review). The most studied of these equations is the nonlocal dimensionless Fisher-KPP equation 
\begin{equation}
\frac{\partial u}{\partial t}= D\frac{\partial^2 u}{\partial x^2} + u\left\{1 - \int_{-\infty}^{\infty} \phi\left(x-y\right) u\left(y,t\right)\,dy \right\}, \label{eqn_NLFKPP0}
\end{equation}
where lengths have been scaled on the nonlocal length scale, and the dimensionless parameter $D$ represents the square of the ratio of the diffusion length scale to the nonlocal length scale. As shown in \cite{BNPR}, this has permanent form travelling wave solutions for all wavespeeds greater than or equal to $2 \sqrt{D}$, with this minimum wavespeed fixed by the behaviour of the solution when $u \ll 1$. The linearisation of (\ref{eqn_NLFKPP0}) when $u$ is small is the same as that of the local Fisher-KPP equation,
\begin{equation}
\frac{\partial u}{\partial t}= D\frac{\partial^2 u}{\partial x^2} + u\left(1-u\right), \label{eqn_FKPP0}
\end{equation}
and the minimum wavespeed exists for the same reason, namely that no strictly positive travelling wave solutions exist for wavespeed less than $2\sqrt{D}$.

Although much of the literature on (\ref{eqn_NLFKPP0}) has focussed on kernels with $\phi(y)>0$ for all $y \in \R$ (see, for example, \cite{JBNL, Gourley2000}), there has also been some interest in kernels with finite support. For example, in \cite{VolpBio, perthame2007,Hamel2014,GVA}, the authors show that for small enough diffusivity such kernels lead to steady solutions with spatial patterning, with either large, narrow spikes (treated as delta functions in \cite{perthame2007}) or spatial patterns with width and height of $O(1)$ as $D \to 0$, which we discuss below. The canonical example of a kernel with compact support is the top hat kernel given below by (\ref{1.8}). Our aim in this series of papers is to develop a thorough understanding of the nature of typical evolutionary dynamics of (\ref{eqn_NLFKPP0}) with the top hat kernel, and the mechanisms involved therein. More generally, systems such as (\ref{eqn_NLFKPP0}) with $D \ll 1$, where the characteristic diffusion length scale is much smaller than the length scale associated with nonlocal interation, have long been known to have solutions with localised spatial patterning, \cite{Gierer1972}, and (\ref{eqn_NLFKPP0}) with $D \ll 1$ is no exception. It is also important to note here that the key features explored in this paper for the top hat kernel are robust to the addition of general piecewise continuous perturbations to this kernel from $L^{\infty}(\mathbb{R}) \cap L^{1}({\mathbb{R}})$, and which are sufficiently small in $L^{1}(\mathbb{R})$. This point will be addressed by the authors  at a later stage in this series of papers.

In order to begin our study, it is convenient to introduce some notation that will be used throughout the paper. 
For any $T>0$ we introduce $\DT\subset\R^2$ by
\begin{align}
    \label{1.1} \DT=\{(x,t)\in\R^2:x\in\R,\,t\in(0,T]\}
\end{align}
with closure $\closedDT$, and also,
\begin{align}
    \label{1.2} \Dinfty=\{(x,t)\in\R^2:x\in\R,\,t\in\R^+\}
\end{align}
with closure $\closedDinfty$.
The Cauchy problem we consider is that concerned with classical solutions  $u:\closedDT\to\R$ to the semilinear, nonlocal, evolution problem, 
\newcommand{\ut}{u_t}
\newcommand{\uxx}{u_{xx}}
\begin{align}
    \label{1.3}     &\ut=D\uxx + u(1-\phi*u),{\on \DT;}
    \\\label{1.4}   &u(x,0)=A g(x)=:u_0(x),\quad\forall\, x\in\R
    \\\label{1.5}   &u(x,t)\to0\text{ as } |x|\to\infty \text{ uniformly for } t\in[0,T].
\end{align}
Here $A>0$, $g\in C(\R) \cap 
 L^\infty(\R)$ and non-negative, $\norm{g}_\infty=1$ whilst $\supp(g)\subseteq[-x_0,x_0]$ $(x_0>0)$. It is also worth noting that the majority of theory developed hereafter will also apply when the initial data has unbounded support, but decays sufficiently rapidly as $|x| \to \infty$ (for example, the Gaussian initial data used in Section~\ref{section2}). Throughout, with $A$, $x_0$ and $g$ prescribed, we will regard a solution to \eqref{1.3}-\eqref{1.5}
as being a solution in the classical sense.
 As discussed above, we further restrict attention to the situation when the nonlocal kernel $\phi:\R\to\R$ has the simple top hat structure, 
\begin{align}
    \label{1.8}\phi(y) & =
	\begin{cases}
	1, & -\frac{1}{2}\leq y\leq \frac{1}{2} \\
	0, & \text{elsewhere},
	\end{cases}
\end{align}
after which, for $(x,t)\in\closedDT$,
\begin{align}
   \label{1.9} (\phi*u)(x,t)=\int_{x-\frac{1}{2}}^{x+\frac{1}{2}}{u(y,t)}dy.
\end{align}

The main focus of the paper will be both the qualitative and quantitative study of the Cauchy problem \eqref{1.3}-\eqref{1.5} with \eqref{1.8} and \eqref{1.9}.
Of particular interest will be the large-$t$ structure of the solution.
For brevity, we will refer to this Cauchy problem as (IBVP) for the rest of the paper.
With this objective in mind, the paper is structured in the following way:\\\\
$\mathbf{In~Section}~\ref{section2}$ we briefly review the {fundamental} questions of uniqueness and global existence for (IBVP), together with some very general basic bounds on the solution. These are readily and rigorously established by applying the results developed in \cite{Hamel2014} to (IBVP).
Then we present detailed illustrative numerical solutions to (IBVP), which enable us to formulate a number of structural and mechanistic conjectures concerning the evolution of the solution to (IBVP) as $t\to\infty$; in particular in relation to the propagation of travelling wavefronts which, below a critical value of $D$ (which is determined explicitly in section 3), leave a stationary spatially periodic state in their wake. Particular attention is paid to the wavelength selection mechanism for this emerging stationary periodic steady state. It is determined, via both numerical computation and large-t asymptotics, that this mechanism changes as the diffusivity $D$ decreases. Specifically, it is shown that for moderately small values of $D$, close to the linear stability margin of the equilibrium state $u=1$, the mechanism approximates to a selection of the most unstable linear wavelength, and operates at a distance to the rear of the wavefronts. However, the mechanism changes significantly for very small values of $D$ into what we refer to as a \emph{hump formation} mechanism, which is controlled by the wavefront and its immediate exponentially small precession, and takes place immediately to its rear, and involves both nonlinear and nonlocal effects. This mechanism is analysed and illuminated in detail via large-t asymptotics.\\\\
$\mathbf{In~Section ~\ref{section4}}$, motivated by section 2 and the references therein, we examine the temporal stability of the two equilibrium states $u=0$ and $u=1$ to the nonlocal  PDE \eqref{1.3} with \eqref{1.9} in detail (via linearization, which is underpinned by classical rigorous linearization theorems away from the marginal stability case), and how this, and its consequences, may relate to (IBVP). In this way, we briefly confirm that the unreacted equilibrium state $u=0$ is temporally unstable at all $D>0$. The linearized analysis at this equilibrium state also confirms that the solution to (IBVP) should develop wavefronts at locations $x \sim \pm2\sqrt{D}t + O(\log{t})$ as $t\to \infty$. This accords with the spreading speeds established in \cite{Hamel2014} (Theorem 1.5), as well as the numerical solutions considered in section 2. We then provide a detailed linearized analysis of the fully reacted equilibrium state $u=1$, and observe that this equilibrium state bifurcates from temporally stable to temporally unstable as $D$ decreases through the critical value $\Delta_1$, with the value $\Delta_1\approx 0.00297$ determined exactly, and in excellent agreement with the numerical solutions presented in section 2. In addition, we present a detailed analysis of the nature of this temporal instability when $0<D<\Delta_1$, which indicates and supports the observations in section 2 that the solution to (IBVP) will develop into a stationary spatially periodic steady state at the rear of the propagating wavefronts. We determine the wavelength of the harmonic Fourier mode which has maximum temporal growth rate and confirm that, \emph{at least at values of $D$ close to the marginal stability value $\Delta_1$}, the selected wavelength of the steady periodic state which emerges in (IBVP) is very close to the wavelength of this Fourier mode with maximum linear growth rate. The analysis in this section affords us the opportunity to propose two fundamental conjectures concerning the large-$t$ asymptotic structure of the solution to (IBVP), which we label as (P1) and (P2). These conjectures motivate the direction of each of the remaining sections of the paper.  \\\\
$\mathbf{In~Section ~\ref{section5}}$, motivated by (P2), we consider in detail the existence of spatially periodic steady states to the nonlocal Fisher-KPP equation featuring in (IBVP), with particular attention devoted to the most interesting situation when $D$ is very small. This section gives a detailed unfolding, in the case of the top hat kernel, of the generic existence result established in \cite{Hamel2014} (Theorem 1.1) regarding spatially periodic steady states. With $\lambda$ representing \emph{fundamental wavelength}, we establish, using local bifurcation/weakly nonlinear theory that a unique (up to spatial translation) nontrivial spatial periodic steady state exists at each point $(\lambda,D) \in \bigcup_{i=1}^{\infty}{\Omega_i}$, where the family $\Omega_i$ are bounded, open, simply connected and pairwise disjoint subsets of the first quadrant of the $(\lambda,D)$ plane, and are constructed explicitly. The periodic steady states are created by a family of  steady state pitchfork bifurcations, from the equilibrium state $u=1$, as the boundary of each subdomain is crossed into its interior. Standard weakly nonlinear approximations to these periodic steady states are recorded for points close to the bifurcation boundary, and these bifurcation curves are path followed numerically to generate the complete bifurcation surface above each subdomain $\Omega_i$. In addition we prove that each of the periodic steady states is strictly positive, and represents an oscillation about $u=1$. Of particular interest is the nature of these periodic steady states as $D\to 0^+$. In this limit we have developed a detailed and intricate theory, via the method of matched asymptotic expansions, to asymptotically approximate the periodic steady states. This enables us to explicitly construct detailed  asymptotic approximations to the periodic steady states as $D\to 0^+$, which are spatially uniform over their wavelength. A key observation is that the structure of the periodic steady states develops into a distinct form of localised hump regions where $u=O(1)$ separated by dead regions where $u$ is exponentially small in $D$, as $D\to 0$. The asymptotic constructions are shown to be in excellent agreement with numerically determined approximations.\\\\
$\mathbf{In~Section ~\ref{section6}}$, in relation to the nature of the propagating fronts identified in the large-$t$ structure of (IBVP), and referred to in (P2), we consider the possibility that the nonlocal Fisher-KPP equation featuring in (IBVP) can support nontrivial, positive, \emph{propagating}, spatially periodic travelling waves which bifurcate from the equilibrium state $u=1$. It is straightforward for us to establish that no such local bifurcations take place, and as such, a possible role of spatially periodic nontrivial travelling  weaves in (IBVP) can be ruled out.\\\\
$\mathbf{In~Section ~\ref{section7}}$ we consider the role in (IBVP) of non-negative travelling waves which have steady profile and represent a transition from the unreacted equilibrium state ahead to the fully reacted equilibrium state to the rear (which we refer to a \emph{transitional permanent form travelling wave} abbreviated to (TPTW) throughout). We establish that such a (TPTW) exists at each $D>0$ when the propagation $v$ speed has $v\ge2\sqrt{D}$, which is in accord with the wavefront evolution speed in (IBVP), as recorded in sections 2 and 3, and in relation to the reported spreading speed for the similar initial value problem reported in \cite{Hamel2014} (Theorem 1.5). The remaining analysis establishes some further elementary properties, and then examines in detail the behaviour to the rear of a (TPTW) with particular emphasis on monotone decay, oscillatory/monotone decay and oscillatory decay. This theory develops, expands and provides explicit critical transition points, for the case of the top hat kernel, in relation to the earlier general theory in \cite{FZ}. \\\\
Finally, in Section \ref{section8}, we bring all of the subsequent results together, with emphasis on their bearing on (IBVP).

\section{General Setting and Numerical Solutions for (IBVP)}\label{section2}

 To begin this section we briefly recall fundamental preliminary results concerning (IBVP), which can be obtained directly from application of the theory developed in \cite{Hamel2014}. 
We first observe, with key details established in \cite{Hamel2014} (Theorem 1.3), that (IBVP) has a unique, global (on $\overline{D}_{\infty}$) solution, which, for each $T>0$, depends continuously on intial data, throughout $\overline{D}_T$. Moreover, $||u(\cdot,t)||_{\infty}$ is uniformly bounded on $[0,\infty)$ by a constant depending only upon $A$ and $D$ (although, for given $A\ge0$ we record from \cite{Hamel2014} (Remark 1.4) that this constant blows up as $D\to 0^+$). In addition, it is trivially established, via the strong maximum principle and comparison theorem, that, for any $t_0$ sufficiently large, 
\begin{align}
    \label{2.3} 0<u(x,t)\le (t+t_0)^{-\frac{1}{2}}e^{(t+t_0)}e^{-x^2/4D(t+t_0)}  \quad \forall\,  (x,t)\in D_T,
\end{align}
which provides useful information in later sections. However, we may observe immediately from \eqref{2.3} that should an identifiable wavefront location and structure develop in the solution to (IBVP) as $t\to \infty$, say at location $|x|\sim S(t)$ (beyond which $u$ is exponentially small as $t\to \infty$), then,
\begin{align}
    \label{2.9} 0<S(t)\leq 2\sqrt{D}t-\frac{1}{2}\sqrt{D}\log t+ O(1),
\end{align}
as $t\to\infty$. We should remark here that Bouin, Henderson and Ryzhik \cite{BHR} establish an asymptotic estimate giving wavefront location at large-$t$ for the solution to the nonlocal Fisher-KPP equation, when the initial data is front-like and localized to the (wlog) left half-line, and in the case of the top hat kernel, this gives a form which satisfies the above inequality, with the coefficient of the $\log t$ term being $-\frac{3}{2}$, as in the Bramson correction for the same evolution problem with the classical local Fisher-KPP equation. We, tentatively, anticipate that in the present situation, $S(t)$ will be similarly asymptotic to this form as $t\to \infty$.

In the remainder of this section we begin our principal study of the detailed qualitative and quantitative features of the solution to (IBVP), and in particular how these properties respond to decreasing the diffusion parameter $D$. We will see that significant changes in structure and mechanism occur, particularly as $D$ decreases from moderately small ($\sim10^{-3}$) to extremely small ($\sim 10^{-6}$). With this in mind, we develop a numerical scheme to approximate solutions to (IBVP), and use this to investigate the qualitative and quantitative structure of solutions to (IBVP), with particular attention to the structure of the solution as $t\to \infty$. We note that similar evolutionary computations have been made and presented in \cite{Nad}, concerning a stability question relating to permanent form travelling wave structures connecting two unstable states of equation (\ref{eqn_NLFKPP0}). Here our aim is to form a preliminary approach to elucidating the qualitative and quantitative properties of the evolution problem (IBVP).
For the first set of numerical solutions that we present, we take $x_0=\frac{1}{2}$ and $g:\R\to\R$ as 
\begin{align}
    \label{3.1}g(x) & =
	\begin{cases}
	(1-2x)^2(1+2x)^2, & |x|\leq \frac{1}{2} \\
	0, & |x|>\frac{1}{2}.
	\end{cases}
\end{align}
We discretise $u$ on a uniform spatial grid of $N$ points, truncated to $0 \leq x \leq L$, approximating the second derivative using central finite differences. The convolution term is evaluated using the trapezium rule, in other words assuming a linear variation of $u$ between grid points, dealing carefully with cases where the edge of the support of the top hat kernel lies between grid points.
We also take into account the symmetry of the solution about $x=0$ and assume that $u=0$ for $x>L$. In the simulations discussed below, we take $L = 10$ and $N = 1000$.
Timestepping is done using the midpoint method (second order accurate), with time step chosen adaptively so that the maximum change in $u$ at each step is below $10^{-2}$.

\begin{figure}
\begin{center}
\includegraphics[width=\textwidth]{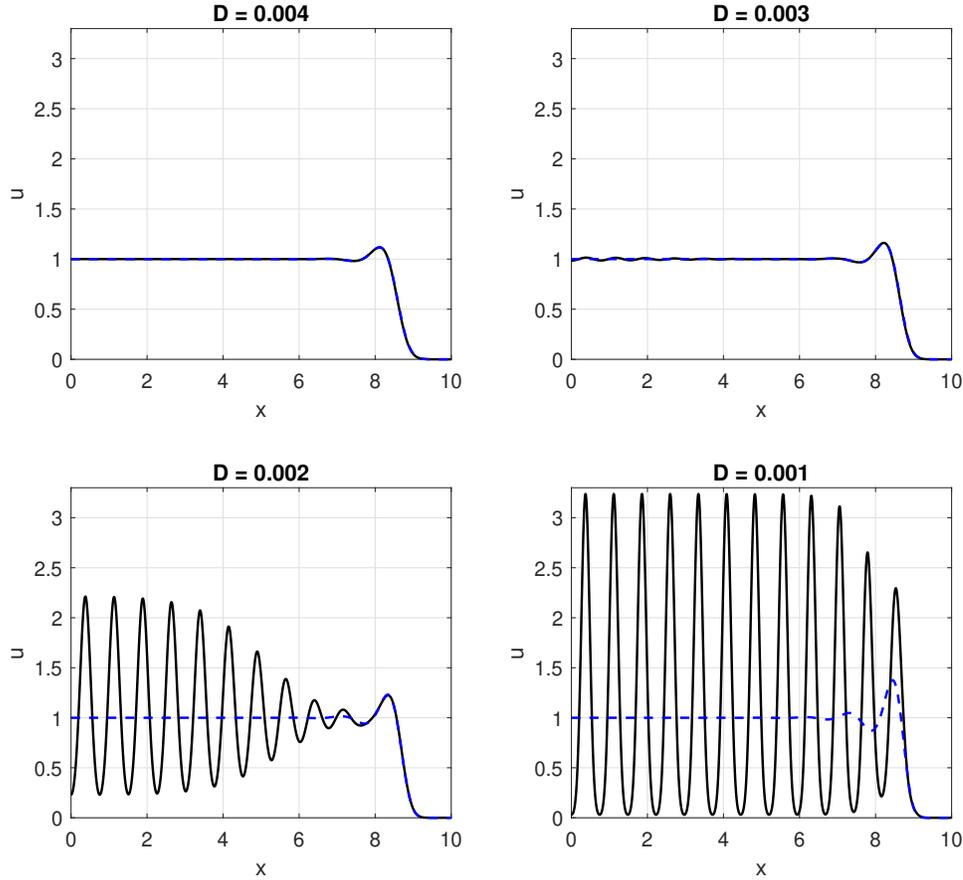}
\caption{The numerical solution of (IBVP) for various values of $D$ with $A =0.01$ (solid black line), along with the minimum speed travelling wave (broken blue line).}\label{fig_TW}
\end{center}
\end{figure}

\begin{figure}
\begin{center}
\includegraphics[width=\textwidth]{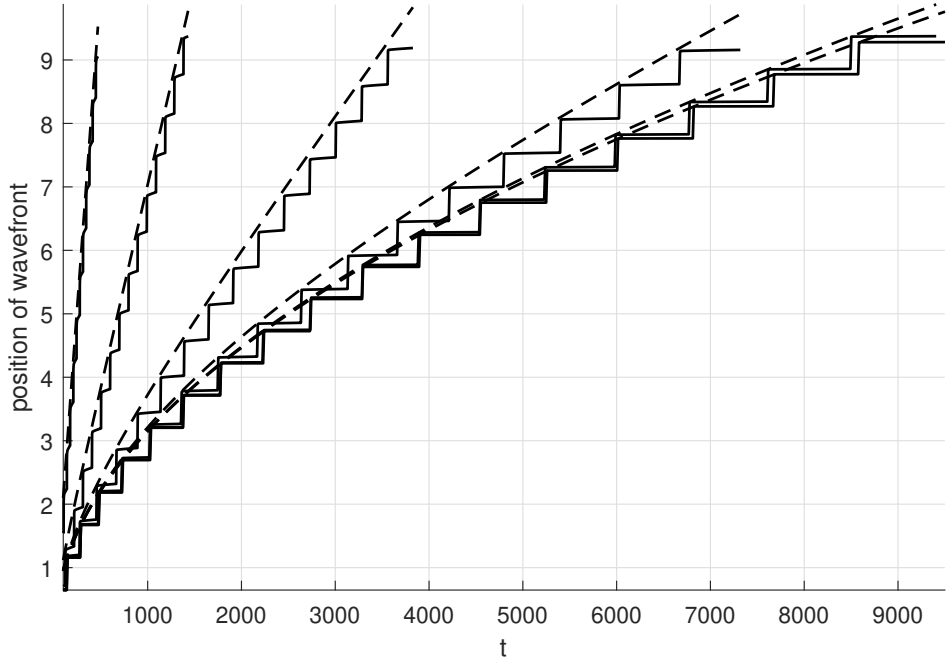}
\caption{The numerically-calculated position of the wavefront for $D = 10^{-4}$, $10^{-5}$, $10^{-6}$, $10^{-7}$, $10^{-8}$ and $10^{-9}$, with $w = 0.1$ and $A=0.01$. The broken line is the function $x_f(t)$, defined in (\ref{eqn_xf}).}\label{fig_TW4}
\end{center}
\end{figure}

\begin{figure}
\begin{center}
\includegraphics[width=\textwidth]{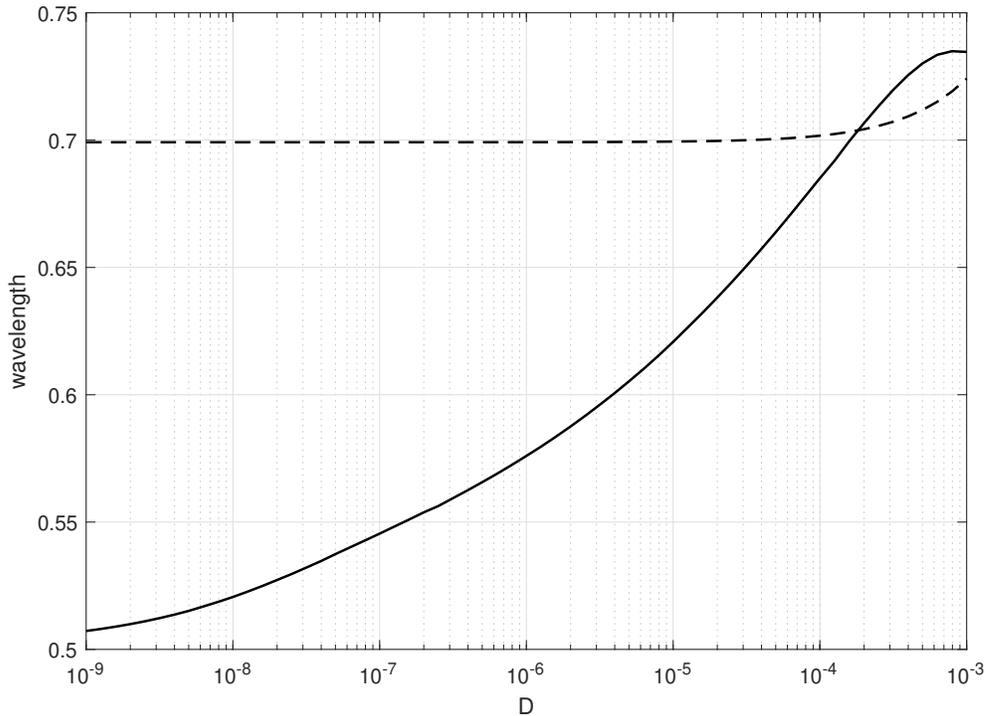}
\caption{The wavelength of the spatially-periodic steady state left behind the wavefront, calculated numerically as a function of $D$. The broken line shows the most unstable wavelength given by the linearized theory.}\label{fig_wavelength}
\end{center}
\end{figure}

Figure~\ref{fig_TW} shows the solution of (IBVP) when $t = 9/2 \sqrt{D}$, which is just before the wavefront reaches the edge of the truncated domain, with $A = 0.01$.
For all values of $A$ that we investigated, indeed for all localised initial inputs of $u$ that we tried, the solution was qualitatively similar to those shown in Figure~\ref{fig_TW}. 
In each case, a wavefront propagates in the positive $x$-direction. 
Also shown is the corresponding minimum speed travelling wave solution (see section \ref{section7} for details of the minimum speed travelling wave solution), calculated numerically using the same finite difference method to set up the discretised equations and `fsolve' in Matlab to solve them. 
The equilibrium state $u=1$ is temporally unstable for $D <\Delta_1  \approx 0.00297$ whilst temporally stable for $D$ larger than this critical value (see section \ref{section4} below).
For $D >\Delta_1$ the minimum speed travelling wave solution emerges, which leaves the equilibrium state $u=1$ in its wake. For $0<D < \Delta_1$ however, a non-propagating stationary, spatially periodic state is left in the wake of the wavefront. For moderately small values of $D$, as shown in Figure~\ref{fig_TW}, the wavelength of this periodic state is close to $0.7$, which is close to the most unstable wavelength in the linear stability theory for the equilibrium state $u=1$, but as we will see, this is only remains so at moderately small values of $D$, with the wavelength selection mechanism and value changing significantly for much smaller values of $D$ (see Section~\ref{section4}, but also Figure~\ref{fig_wavelength} and subsequent discussion below). Movies of the numerical solutions illustrated in Figures~\ref{fig_TW} and~\ref{fig_TW4} can be found \href{https://drive.google.com/drive/folders/1Q8MjbXVh1eRuSO_fxEY-7OWqNLy1MAlB?usp=sharing}{here}.

\begin{figure}
\begin{center}
\includegraphics[width=\textwidth]{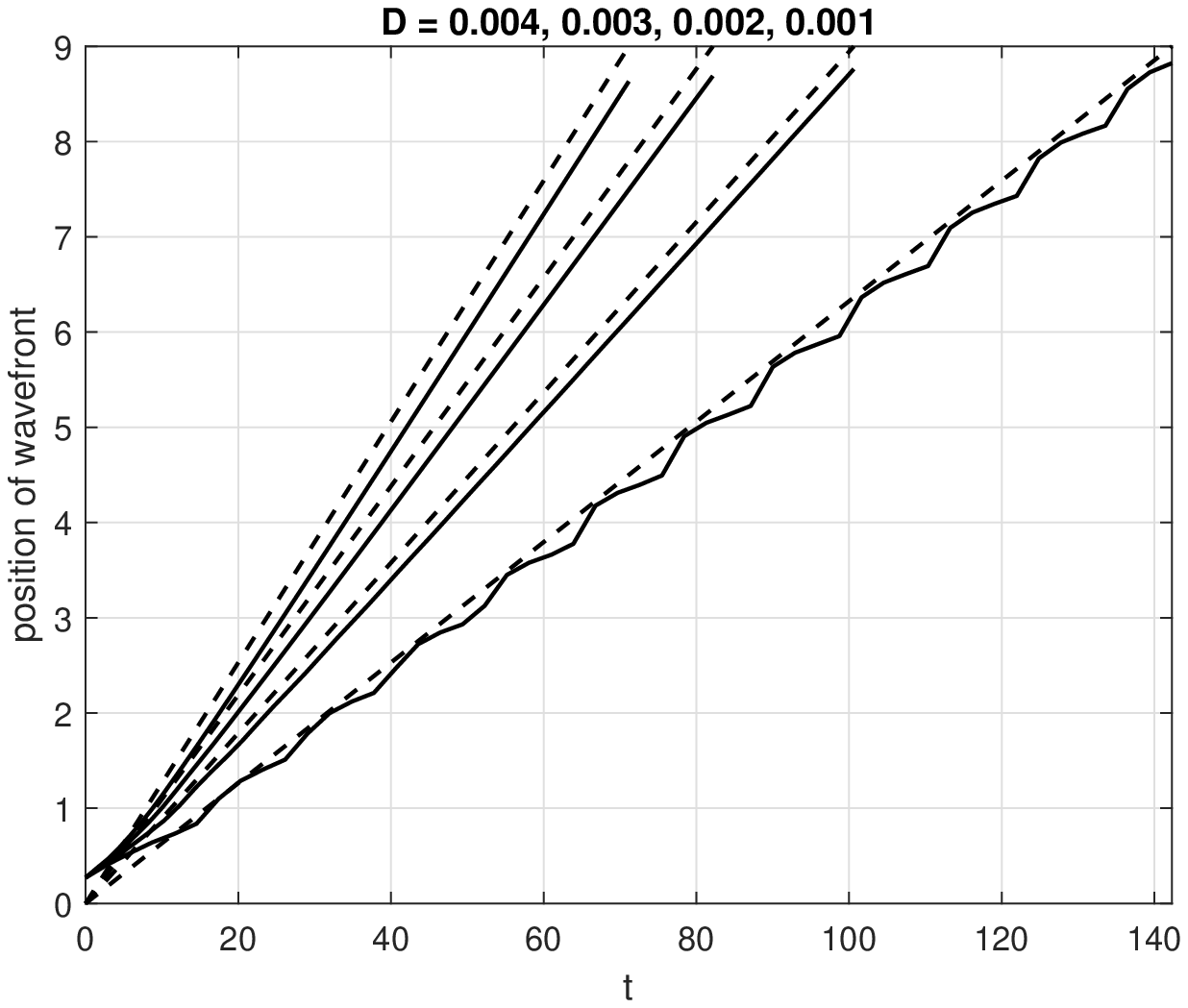}
\caption{The numerically-calculated position of the wavefront for various values of $D$. The broken line has slope $2 \sqrt{D}$, the minimum wavespeed.}\label{fig_TW2}
\end{center}
\end{figure}

The temporal periodicity of the creation of the stationary spatially periodic state behind the wavefront is also reflected in a weak periodic variation in the position of the wavefront (defined to be the largest value of $x$ at which $u=\frac{1}{2}$), but with an average speed of approximately $2\sqrt{D}$, the minimum wavespeed, which is clearly illustrated in  Figure~\ref{fig_TW2} (see sections \ref{section4} and \ref{section7}, regarding the specific notion of minimum wavespeed in the present context).
For the solutions shown in Figure~\ref{fig_TW2} with $D=0.003$ and $D=0.002$, the formation of a stationary periodic state behind the wavefront does not cause oscillations in the position of the wavefront, whereas for $D = 0.001$ the magnitude of the oscillation is large enough that there is a weak effect. Although choosing a smaller value of $u$ at which to define the position of the wavefront would eliminate these traces of the oscillations, as $D \to 0$, the value of $u$ required to do so becomes exponentially small. It therefore seems reasonable to use $u = \frac{1}{2}$ as our qualitative definition of the position of the wavefront, and thereby retain an indication of the fundamentally oscillatory nature of the observable solution in our record of the progress of the travelling wave.

\begin{figure}
\begin{center}
\includegraphics[width=\textwidth]{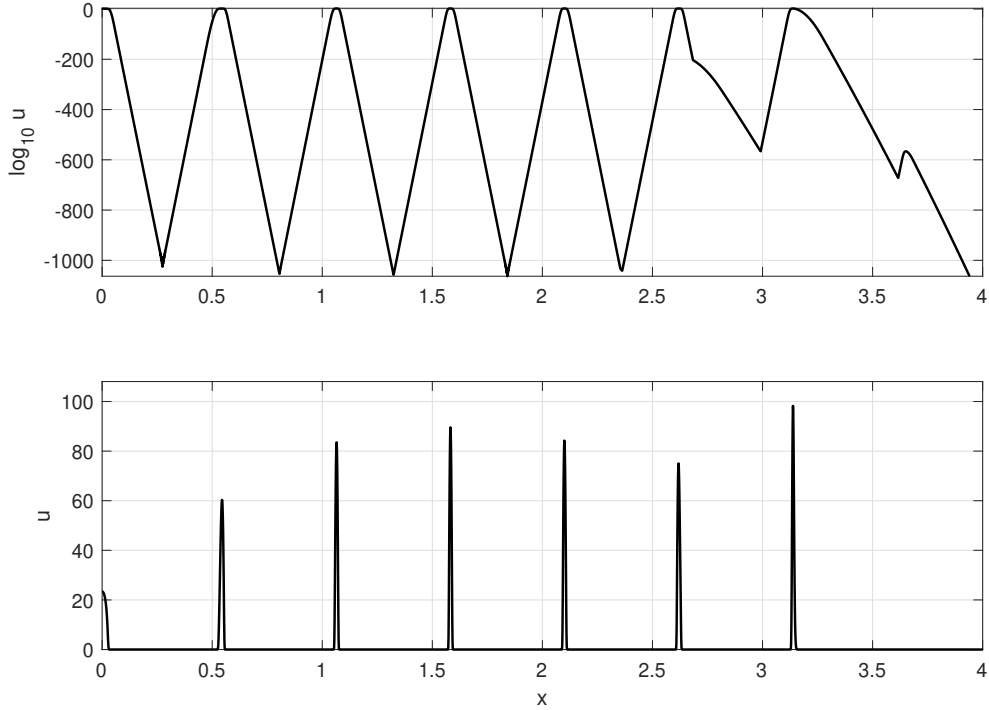}
\caption{The numerical solution of (IBVP) for Gaussian initial data with $A=0.01$ and width $w = 0.04$, and $D = 10^{-8}$, when $t = 6000$. The upper panel shows $\log_{10} u$. New spikes are initiated ahead of the wave at the point where $u$ is close to $10^{-700}$, which can only be captured accurately by solving for $\log u$ instead of $u$.} \label{fig_TW3}
\end{center}
\end{figure}

\begin{figure}
\begin{center}
\includegraphics[width=\textwidth]{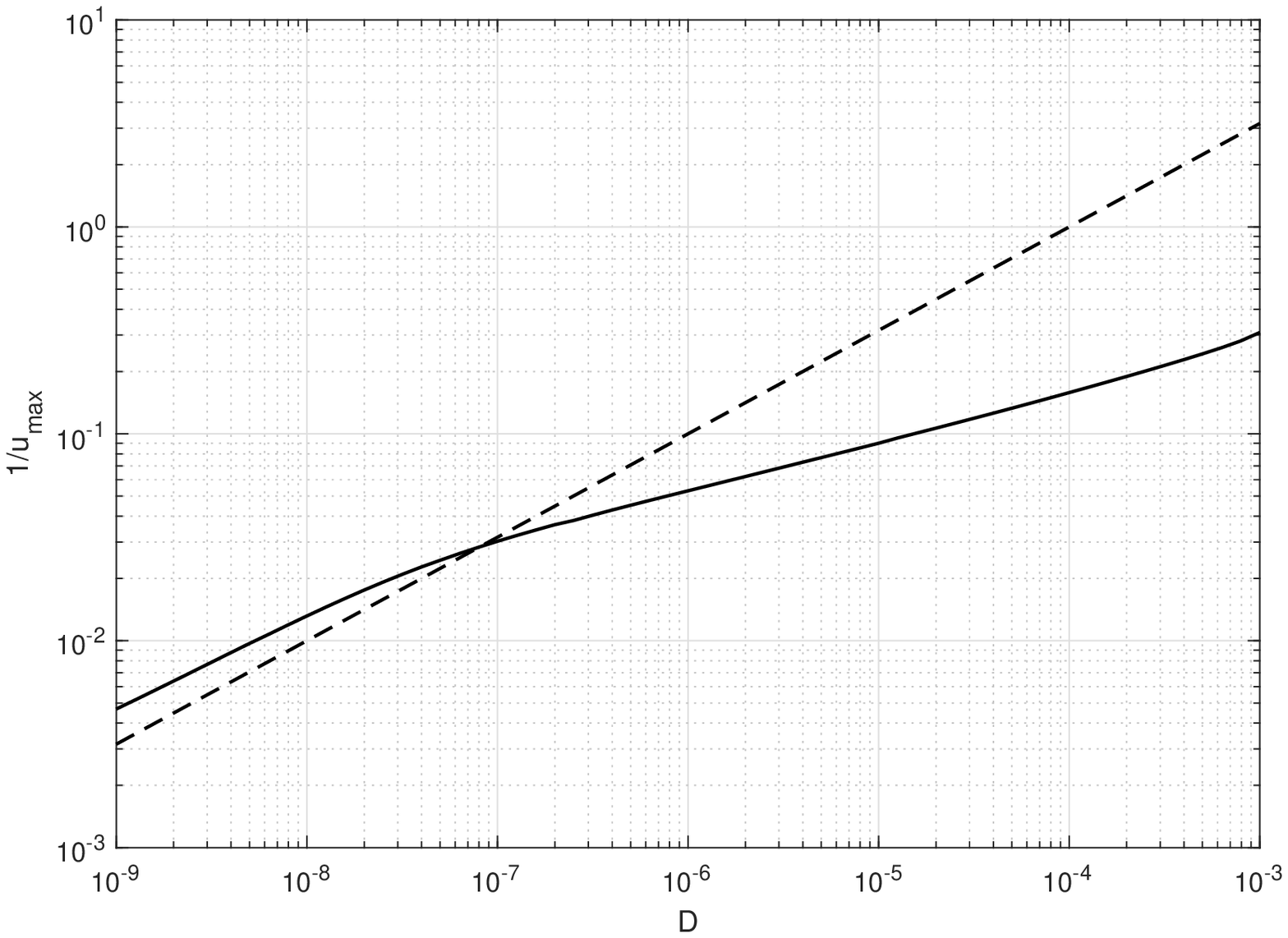}
\caption{The inverse of the height of the spikes behind the wavefront, calculated numerically as a function of $D$. The broken line is 100/$\sqrt{D}$.}\label{fig_umax}
\end{center}
\end{figure}

For smaller values of $D$ than those used in Figures~\ref{fig_TW} and~\ref{fig_TW2}, we find that the creation of the humps, which ultimately form the stationary and spatially periodic state to the rear of the wavefront, initiates, periodically in $t$, just ahead of the wavefront and, as we shall see below, can be related to the dynamics of that part of the solution profile that becomes exponentially-small with distance ahead of the wavefront. It is now this periodic mechanism that selects the final spatial wavelength of the stationary periodic state which forms at the rear of the wavefront. In order to accurately compute this part of the solution, we use a different numerical method, solving instead for $W = \log u$, which satisfies the evolution equation
\begin{equation}
    \frac{\partial W}{\partial t} = D\frac{\partial^2 W}{\partial x^2} + D \left(\frac{\partial W}{\partial x}\right)^2 + 1 - \phi * e^W. \label{eqn_W}
\end{equation}
We also use an FFT to calculate the convolution term (and periodic boundary conditions with periodicity large enough that the effect on the travelling wave dynamics is negligible) and use five point stencils for the derivatives for greater accuracy. In addition, since we cannot take the logarithm of an initial condition with compact support, we instead use a Gaussian of width $w>0$ as the initial condition, so that
\begin{equation}
    g(x) = e^{-x^2/w^2}.
\end{equation}
Figure~\ref{fig_TW3} shows a snapshot of the evolving solution for $D = 10^{-8}$, $A=0.01$ and $w = 0.1$. 
Although the creation of a spatially periodic steady state behind the wavefront is qualitatively similar to that shown for larger values of $D$ in Figure~\ref{fig_TW}, we can now see how the behaviour ahead of the wavefront leads to the creation of the humps, which have unit weight, width of $O(D^{1/2})$ and height of $O(D^{-1/2})$ for $D \ll 1$ (see Section~\ref{sec_5.3.3} and Figure~\ref{fig_umax}). Far ahead of the wavefront, the solution evolves as
\begin{equation}
    u(x,t) \sim \frac{w A}{\sqrt{w^2+4Dt}} \exp\left(t - \frac{x^2}{w^2+4Dt}\right)~~\mbox{for $x \gg t \geq O(1)$.}\label{eqn_ff}
\end{equation}
Note that (\ref{eqn_ff}) complies with the upper bound given by (\ref{2.3}) with $t_0 = w^2/4D$. In a spatial interval of unit width, centred on the leading hump, the term $1-\phi*u$ is small. However, ahead of this interval, it becomes close to one. This has the effect of turning on the term $e^t$ in (\ref{eqn_ff}), which manifests itself in the solution shown in Figure~\ref{fig_TW3} as a rapid change in $\log_{10}u$, ahead of which $u$ grows until a new hump is formed. This is most easily observed in the animation that can be found \href{https://drive.google.com/drive/folders/1Q8MjbXVh1eRuSO_fxEY-7OWqNLy1MAlB?usp=sharing}{here}. This hump formation process, generated by the dynamics of the exponentially small part of the solution ahead of the wavefront, is similar to that studied in \cite{JB1} for a nonlocal reaction-diffusion equation with a different local reaction term. For $D\ll 1$, this hump forms when the logarithm of the far field solution (\ref{eqn_ff}) becomes small, and we can approximate this location as the point where the logarithm is zero, namely $x = x_f(t)$, with
\begin{equation}
    x_f(t) = \sqrt{w^2+4Dt} \sqrt{t + \log(Aw) - \frac{1}{2} \log(w^2+4Dt)}.\label{eqn_xf}
\end{equation}
Figure~\ref{fig_TW4} shows the position of the wavefront along with $x_f(t)$ for various values of $D$. 
As can be seen, $x_f(t)$ gives an excellent prediction of the position of the wavefront, and we observe that it conforms with the wavefront bound given in \eqref{2.9}. We also note that this hump formation process naturally leads to a stationary, spatially periodic solution behind the wavefront with wavelength close to $\frac{1}{2}$. Figure~\ref{fig_wavelength} shows how the wavelength varies with $D$, consistent with this observation. 
Also shown is the most unstable wavelength according to the linearized theory of Section~\ref{section4}, which is close to the wavelength observed for moderately small values of $D$, when the periodic states emerges behind the wavefront, \emph{but not for smaller values of} $D$, where the hump formation mechanism described above takes over and creates this state, and determines its wavelength to be approaching $\frac{1}{2}$ as $D\to 0^+$. Finally, it is anticipated that additional details relating to this hump formation mechanism may be made available by a more detailed development of the large-$t$ asymptotic form \eqref{eqn_ff} via the method of matched asymptotic coordinate expansions (see Leach and Needham \cite{LNB}) or the approximation methods detailed in van Saarloos \cite{Van}. However, for the present paper we regard the above details as sufficient to elucidate this particular mechanism.
Many of the additional  distinct features described above, particularly when the diffusivity is very small, have not been reported in earlier literature, and their illumination is the focus of the following sections.

\section{Equilibrium States, Stability Characteristics and linear evolution}\label{section4}
The nonlocal PDE \eqref{1.3}, with \eqref{1.9}, has two equilibrium states. The \emph{{unreacted}} state with 
\begin{align}
    \label{4.1} u(x,t)=0\quad \forall\,\xt\in\closedDinfty,
\end{align}
and the \emph{{fully reacted}} state 
\begin{align}
    u\xt =1\quad \forall\,\xt\in\closedDinfty.
\end{align}
We have seen in Section \ref{section2} that these equilibrium states play a key role in the large-$t$ structure of the solution to (IBVP).
Of particular significance is the temporal stability of these equilibrium states.
In this section we examine the linearized temporal stability of each of these equilibrium states in detail, and consider the consequences thereof in terms of spatio/temporal evolution.
To this end, we formulate a linearized initial value problem.
We write,
\newcommand{\ubar}{\overline{u}}
\newcommand{\ubart}{\ubar_t}
\newcommand{\ubarxx}{\ubar_{xx}}
\begin{align}
   \label{4.3} u(x, t)=u_{e}+\delta \ubar(x, t),\quad (x, t) \in \closedDinfty,
\end{align}
with $\delta \ll 1$ and $u_e=1$ when considering the fully reacted state or $u_e=0$ when considering the unreacted state.
On substituting from \eqref{4.3} into \eqref{1.3}, and neglecting terms of $O(\delta^2)$ as $\delta\to0$, we obtain a linear evolution equation for $\overline{u}$, namely
\begin{align}
    \label{4.4}\ubart=D\ubarxx +L(\ubar), {\on\Dinfty},
\end{align}
where 
\begin{align}
    \label{4.5}L(\ubar)=
    \begin{dcases}
    \ubar;& u_e=0\\
    -\int_{x-\frac{1}{2}}^{x+\frac{1}{2}}\ubar(y,t)dy;& u_e=1
    \end{dcases}
\end{align}
after using \eqref{1.9}.
The linearized initial value problem is then composed of \eqref{4.4}, with associated initial and far field conditions,
\begin{align}
    \label{4.6} &\ubar (x,0)=\overline{g}(x)\quad\forall\,x\in\R
    \\\label{4.7}&\ubar\xt\to 0 \text{ as } |x|\to\infty \text{ uniformly on }\closedDT \text{ (each }T>0).
\end{align}
Here $\overline{g}\in C^1(\R) \cap 
 L^\infty(\R)$ and non-negative, $\norm{\overline{g}}_\infty=1$ whilst $\supp(\overline{g})\subseteq[-x_0,x_0]$ $(x_0>0)$. This problem will be referred to as \livp0 when $u_e=0$, and \livp1 when $u_e=1$. We remark that aspects of the following linearized theory have been discussed in \cite{VolpBio} in relation to equilibrium state temporal stability. Here we develop this further in terms of detailing the full qualitative structure associated with the linearized evolution problems, which provides additional and significant information towards the analysis of (IBVP).
We now consider \livp0 and \livp1 in turn.

\subsection{Analysis of (LIVP)\texorpdfstring{\(_0 \)}.}

We first seek elementary solutions to \eqref{4.4} and \eqref{4.5} in the form,
\begin{align}
    \label{4.8} \ubar\xt=e^{ikx-w t},\quad \xt\in\closedDinfty
\end{align}
with $k\in\R$ and $w\in\C$. 
On substitution from \eqref{4.8} in \eqref{4.4} and \eqref{4.5} we obtain the dispersion relation 
\begin{align}
    \label{4.9} w=w_0(k)=Dk^2-1\quad\forall\,k\in\R
\end{align}
and so, as expected, \livp0 is \emph{nondispersive} ($w_0(k)\in\R$ $\forall$  $k\in\R$), and moreover, for each $D>0$,
\begin{align}
    \label{4.10}w_0(k)<0\quad\forall\,k\in\left(\frac{-1}{\sqrt {D}},\frac{1}{\sqrt{D}}\right)
\end{align}
with
\begin{align}
    \label{4.11}w_0(0)=-1=\inf_{k\in\R}w_0(k).
\end{align}
We immediately conclude that the equilibrium state $u_e=0$ is temporally unstable, according to the linearized theory, at each $D>0$
(this conclusion can readily extended to apply to the fully nonlinear and nonlocal PDE \eqref{1.3} with \eqref{1.9}, via an application of the parabolic comparison theorem to the operator $N(w):=w_t-Dw_{xx}-w$ on $\DT$; for brevity we omit the details). The Fourier integral theorem allows us to write down the solution to {\livp0} and this can then be estimated, via steepest descents, to obtain, 
\begin{align}
    \label{4.14}\ubar \xt \sim \frac{\sqrt{\pi}}{\sqrt{D}t^{\frac{1}{2}}}\widehat{g} (0)e^te^{-x^2/4Dt}
\end{align}
as $t\to \infty$ uniformly for $x\in\R$, with
\begin{align}
    \notag \widehat{g}(0)=\frac{1}{2\pi}\int_{-x_0}^{x_0}\overline{g}(s)ds \quad(>0).
\end{align}
We observe from \eqref{4.14} that, as $t\to\infty$, there are two symmetric `wavefronts' where 
\begin{align}
    \label{4.15} |x|\sim 2\sqrt{D}t\text{ as } t\to\infty,
\end{align}
and behind the wavefronts $\ubar$ is growing exponentially in $t$, whilst ahead of the wavefronts $\ubar$ is decaying exponentially in $t$. 
This structure is consistent with \eqref{2.9}, and the numerical solutions to (IBVP) in section \ref{section2}. It is also consistent with the rigorous result in \cite{BHR} for the corresponding evolution problem when the initial data is front-like and localized to the left half-line.

\subsection{Analysis of (LIVP)\texorpdfstring{\(_1 \)}.}

We seek elementary solutions to \eqref{4.4} and \eqref{4.5} in the form of \eqref{4.8}, with again $k\in\R$ and $w\in \C$.
This now leads directly to the dispersion relation 
\begin{align}
    \label{4.16} w=w_1(k)=Dk^2+\frac{2}{k}\sin\frac{1}{2}k\quad\forall\,k\in\R
\end{align}
and we observe that $w_1(k)$ is an even function of $k$.
Moreover $w_1(k)\in\R$ for all $k\in\R$ and so \livp1 is \emph{nondispersive}.
In further analysing \eqref{4.16}, it is convenient to introduce the function $\Delta:\R^+\to\R$ such that 
\begin{align}
    \label{4.17} \Delta(X) =-\frac{2}{X^3}\sin\frac{1}{2}X\quad\forall\,X\in\R^+.
\end{align}
The zeros of $\Delta(X)$ are at
\begin{align}
    \label{4.18} X=2n\pi,\quad n\in\N,
\end{align}
whilst the turning points are at 
\begin{align}
    \label{4.19} X=\delta_n,\quad n\in\N,
\end{align}
with $2n\pi<\delta_n<2(n+1)\pi$, and 
\begin{align}
    \label{4.20}\delta_n\sim(2n+1)\pi \text{ as } n\to\infty.
\end{align}
Furthermore, $\delta_n$ is a local maximum when $n$ is odd and a local minimum when $n$ is even. 
At each local maximum point we write,
\begin{align}
    \label{4.21} \Delta_r =\Delta(\delta_{2r-1}),\quad r=1,2,\dots.
\end{align}
We observe that,
\begin{align}
    \label{4.22} 0<\Delta_{r+1}<\Delta_r,\quad r=1,2,\dots
\end{align}
and 
\begin{align}
    \label{4.23}\Delta_r\sim \frac{2}{\pi^3 (4r-1)^3}\text{ as }r\to\infty,
\end{align}
whilst a straightforward numerical calculation gives ${\Delta}_1\approx0.00297$. 

We can now readily interpret the dispersion relation \eqref{4.16}. We first observe that for 
\begin{align}
    \label{4.24}D>\Delta_1
\end{align}
then $w_1(k)>0$ for all $k\in\R$ and so the equilibrium state $u_e=1$ is temporally asymptotically stable.
However, for
\begin{align}
    \label{4.25} 0<D<\Delta_1
\end{align}
then
\begin{align}
    \label{4.26}\inf_{k\geq0}w_1(k)=w_1(k_m)=k_m^2(D-\Delta(k_m))<0,
\end{align}
and so the equilibrium state $u_e=1$ is now temporally unstable. 
We note that $k=k_m$ is uniquely determined as the smallest positive root of the transcendental equation,
\begin{align}
    \label{4.27}w_1'(k) \equiv 2k(D-\Delta(k)) - k^2\Delta'(k) \equiv (2Dk^3-2\sin\frac{1}{2}k +k\cos\frac{1}{2}k)k^{-2}=0
\end{align}
and it is straightforward to establish from this that,
\begin{align}
\label{4.28}k_m\to\begin{cases}
\delta_1\text{ as }&D\to\Delta_1^-\\
k_0\text{ as } &D\to 0^+
\end{cases}
\end{align}
where $k=k_0$ is the smallest positive root of the equation 
\begin{align}
    \label{2.49} \tan\frac{1}{2}k=\frac{1}{2}k
\end{align}
so that $2\pi<k_0<3\pi$, and
\begin{align}
    \label{4.30}w_1(k_m)\to\frac{2}{k_0}\sin\frac{1}{2}k_0\text{ as }D\to0^+.
\end{align}
We note that in both cases,
\begin{align}
    \notag w_1(k)\sim Dk^2\text{ as }|k|\to\infty
\end{align}
and so \livp1 is well-posed. Figure~\ref{fig_wavelength} shows the wavelength of the most unstable mode as a function of $D$.

The solution to \livp1 is readily obtained as
\begin{align}
    \label{4.31} \ubar\xt=\int_{-\infty}^\infty \widehat{g}(k)e^{-w_1(k)t}e^{ikx}dk
\end{align}
for $\xt\in\closedDinfty$, with $\widehat{g}$ being the Fourier transform of $\overline{g}$.
We obtain from \eqref{4.31}, via Laplace's method, that, 
\begin{align}
    \label{4.32} \ubar\xt\sim\frac{\sqrt{\pi}2^{3/2}|\widehat{g}(k_m)|}{(w_1''(k_m))^{\frac{1}{2}}t^{\frac{1}{2}}}e^{-w_1(k_m)t} 
    \cos{(k_m x+\arg( \widehat{g}(k_m))}
\end{align}
for $|x|=O(1)$ as $t\to\infty$.
It should be noted that further regions in the large-$t$ structure of $\ubar$ are required when $|x|=O(t)$ as $t\to\infty$, which gives the transition into the spatially exponentially decaying far field for $|x|\gg O(t)$.
We observe, from \eqref{4.32}, that when $0<D<\Delta_1$, the solution to \livp1  evolves into an exponentially growing harmonic periodic state with spatial wave number $k_m$, which depends upon $D$, and temporal exponential growth rate $|w_1(k_m)|$.
This periodic state is stationary, and evolves when $|x|=O(1)$ as $t\to\infty$, behind transition into the far field when $|x|=O(t)$ as $t\to\infty$.\\

In relation to (IBVP), the above analyses indicate that when $t$ is large, two symmetric permanent form wavefronts propagate to left and right, with asymptotic propagation speeds of $\pm2\sqrt{D}$. However, \emph {in the region to the rear of the two wavefronts at} $|x|\sim 2\sqrt{D}t$, the large-$t$ structure to \livp1 indicates that, as $t\to\infty$ in (IBVP), and specifically when $|x|=O(1)$ as $t\to\infty$, \emph{there are two possible steady spatial structures which emerge as a consequence of wavefront passage}, depending upon $D$. 
With $u:\closedDinfty\to\R$ being the solution to (IBVP), these possibilities are:\\
\begin{itemize}
    \item[(P1)] when $D>\Delta_1$, then $u\xt\to1$ as $t\to\infty$, uniformly with $|x|=O(1)$\\
    \item[(P2)] when $0<D<\Delta_1$, then $u\xt\to P(x)$ as $t\to\infty$, uniformly with $|x|=O(1)$. 
    Here $P:\R\to\R$ is a steady periodic solution to equation \eqref{1.3} with \eqref{1.9}, which is positive, has $1\in\Im(P)$ and has wavelength which \\
    (i) is close to $2\pi/k_m(D)$ \emph{when $D$ is close to $\Delta_1$}, and is selected well behind the emerging wavefronts by a \emph{maximum linear growth rate mechanism} on the now weakly linearly unstable equilibrium state $u=1$,\\
    and\\
    (ii) moves away from this value as $D$ decreases, decreasing and drifting towards the minimum wavelength of $1/2$ (see Figure~\ref{fig_wavelength} and section 4) as $D\to 0^+$, with \emph{this wavelength now selected via the hump formation mechanism}, which operates to the immediate rear of the wavefronts (see section 2).\\
\end{itemize}
We recall that both (P1) and (P2) are supported by the numerical solutions to (IBVP) presented in section \ref{section2} in conjunction with the detailed linearized theory of this section.

As a consequence of the theory in this section, and in particular to further enable the development of theory to support and extend the conjecture (P2), the next natural key step is to investigate the existence and structural nature of the positive periodic steady solutions to the nonlocal equation \eqref{1.3} with \eqref{1.9}, which are conjectured to emerge in the large-$t$ development of the solution to (IBVP).
The starting point for this study begins by investigating the emergence of periodic steady solutions via steady state bifurcations from the equilibrium solution $u_e=1$. Thereafter, the detailed structure of these fully nonlinear and nonlocal periodic solutions is developed in detail for significant asymptotic limits. In particular, a fully developed theory is obtained in the small $D$ limit, when the nonlocal length scale is much larger than the diffusion length scale, leading to spatially periodic states consisting of separated, periodically distributed, localised humps, characterised by nonlocal effects which are regulated by weak diffusion.

\section{Positive Periodic Steady States}\label{section5}

We begin by considering the steady state form of the nonlocal equation \eqref{1.3} with \eqref{1.9}, and in particular we seek to identify steady state bifurcations from the equilibrium state $u_e=1$, which give rise to periodic steady states. For the general form of the nonlocal Fisher-KPP equation, under quite general conditions on the kernel, generic considerations concerning local bifurcations to small amplitude bounded steady states on the real line, as the associated equilibrium state loses its temporal stability, have been developed in \cite{Fay} and a detailed consideration of the spectral theory on the real line of the linear operators thereby encountered has been made by Volpert and Vougalter \cite{VV}. By restricting attention to bifurcations to \emph{periodic steady states}, we effectively restore compactness to the operators, and can then develop a standard weakly nonlinear bifurcation theory. Indeed, under the periodic restriction, these operator theoretic results conform with the basic weakly nonlinear theory outlined below for the current specific situation, which we then develop in detail onto the global bifurcation branch.
Any steady state to the nonlocal equation \eqref{1.3} with \eqref{1.9}, in the present context, is a function $F:\R\to\R$, with $F\in C^{2}(\R)\cap L^\infty(\R)$, and such that,
\begin{align}
    \label{5.1} DF_{xx}+F\left(1-\int_{x-\frac{1}{2}}^{x+\frac{1}{2}}F(y)dy\right)=0,\quad x\in\R.
\end{align}
We restrict attention to considering the existence of positive periodic steady states, which oscillate about the equilibrium state $u_e=1$.
Specifically, we will fix $D>0$, and consider the bifurcation to periodic steady states from $u_e=1$, with fundamental wavelength $\lambda$, as the bifurcation parameter.
As a preliminary we first observe, via bootstrapping in \eqref{5.1}, that when $F\in C^2(\R) \cap L^\infty (\R)$ is a steady state, then in fact, $F\in C^\infty(\R)$.
Moreover (see for example \cite{Hale1977}) this can be improved to 
\begin{align}
    \label{5.2} F\in C^\omega(\R) \cap L^\infty(\R).
\end{align}

Now, let $F=F_p(x,\lambda,D)$ be a positive, periodic, steady state at diffusivity $D$ and with fundamental wavelength $\lambda>0$. 
We define
\begin{align}
    \label{5.3} \alpha \equiv \max_{x\in[0,\lambda)}F_p(x,\lambda,D)-\min_{x\in[0,\lambda)}F_p(x,\lambda,D).
\end{align}
the peak-to-trough magnitude of this periodic steady state.
In general, we anticipate that $\alpha=\alpha(\lambda,D)$.
For fixed $D>0$, we now suppose that a bifurcation to periodic steady states occurs from $u_e=1$ as $\lambda$ passes through $\lambda=\lambda_b$ $(>0)$.
To determine the possible values of $\lambda_b$, we consider those values $\lambda=\lambda_b$ when the linearized form of \eqref{5.1} has solution 
\begin{align}
    \label{5.4} F(x)=1+\alpha \cos(\frac{2\pi}{\lambda_b}x+\phi),\quad x\in\R,
\end{align}
as $\alpha\to0^+$, with $\phi$ being a constant phase.
On substitution from \eqref{5.4} into \eqref{5.1}, and allowing $\alpha\to0^+$, a non-trivial solution requires that $\lambda_b$ satisfies the transcendental equation, with fixed $D>0$,
\begin{equation}
    D =-\frac{\lambda_b^3}{4\pi^3}\sin\frac{\pi}{\lambda_b}
    =\Delta\left(\frac{2\pi}{\lambda_b}\right)\label{5.5}
\end{equation}
with $\Delta(\cdot)$ as introduced in \eqref{4.17}.
We can now interpret the root structure of \eqref{5.5}, using \eqref{4.18}-\eqref{4.24}.
In particular, for each $r=1,2,\dots$, then, with 
\begin{align}
    \label{5.6}\Delta_{r+1}<D<\Delta_r,
\end{align}
equation \eqref{5.5} has exactly $2r$ positive roots, which we label as $\lambda_i^\pm(D)$ for \
{$i=1,2,\dots,r$}, with 
\begin{align}
    \label{5.7} \frac{1}{2i}<\lambda_i^-(D)<\frac{2\pi}{\delta_{2i-1}}<\lambda_i^+(D)<\frac{1}{(2i-1)}
\end{align}
for each $i=1,\dots,r$.
We note that, for fixed $k\in\N$, then $\lambda^\pm_k(D)$ exist and are continuous for $D\in[0,\Delta_k)$.
Moreover, the following properties are readily established:
\begin{align}
    \label{5.8}& \lambda_k^\pm (D)\to\frac{2\pi}{\delta_{2k-1}}\text{ as } D\to\Delta_k^- ,
    \\\label{5.9}& \lambda_k^{+(-)}(D)\text{ is monotone decreasing (increasing) with } D\in[0,\Delta_k),
    \\\label{5.10}&\lambda_k^+(D)\to\frac{1}{(2k-1)}^-\text{ as }D\to0^+,
    \\\label{5.11}&\lambda_k^-(D)\to\frac{1}{2k}^+\text{ as }D\to0^+.
\end{align}
Now, fix $r\in\N$, and fix $\Delta_{r+1}<D<\Delta_r$.
At this fixed $D$, there are thus $2r$ local bifurcation points to periodic steady states, namely at critical wavelengths 
\begin{align}
    \label{5.12} \lambda=\lambda_k^\pm(D),\quad  k=1,2,\dots,r.
\end{align}
Further lengthy, but straightforward calculations (standard weakly nonlinear asymptotics as the bifurcation point is approached) establish (which can be made rigorous in a standard way  within the classical local bifurcaton theory, as developed in, for example, Kuznetsov \cite{KUZ})  that each bifurcation is a steady state pitchfork bifurcation from $u_e=1$, being subcritical at $\lambda=\lambda_k^+(D)$ and supercritical at $\lambda=\lambda_k^-(D)$ $(k=1,2,\dots,r)$.
Additional investigation (numerical continuation of the weakly nonlinear bifurcation branches) reveals that, for each $k$, the two bifurcated branches at $\lambda_k^\pm(D)$ are connected in the $(\lambda,\alpha)$ bifurcation diagram.
We can represent this curve as $\alpha=\alpha(\lambda;D)$ for $\lambda_k^-(D)\leq \lambda \leq \lambda_k^+(D)$ (for each $k=1,2,\dots,r$).
Close to the bifurcation points $\lambda=\lambda_k^\pm(D)$ a very standard weakly nonlinear local bifurcation theory is readily developed (details are omitted for brevity as they follow a now standard routine) to establish that, as should be expected for a nondegenerate bifurcation,
\begin{align}
    \label{5.13}\alpha(\lambda,D)=C_k^\pm(D)\abs{\lambda-\lambda_k^\pm(D)}^{\frac{1}{2}}+O(\abs{\lambda-\lambda_k^\pm(D)})
\end{align}
as $\lambda\to(\lambda_k^+(D))^-$ or as $\lambda\to(\lambda_k^-(D))^+$ respectively. 
Here $C_k^\pm(D)$ is a positive constant.
The bifurcated periodic steady states have the form
\begin{align}
    \label{5.14} u=F_p(x,\lambda,D)=1+\alpha(\lambda,D)\cos(\frac{2\pi}{\lambda}x+\phi)+O(\alpha^2(\lambda,D))
\end{align}
as $\lambda\to\lambda_k^\pm(D)$ respectively, with $\alpha(\lambda,D)$ as in \eqref{5.13}, and $\phi\in[0,2\pi)$ being an arbitrary phase (representing translational invariance of \eqref{5.1} in $x$).
In addition, numerical investigation confirms that, with $D$ fixed, then $\alpha(\lambda,D)$ has a single stationary point for $\lambda\in[\lambda_k^-(D),\lambda_k^+(D)]$, which is a maximum point.

\begin{figure}
\begin{center}
\includegraphics[width=0.8\textwidth]{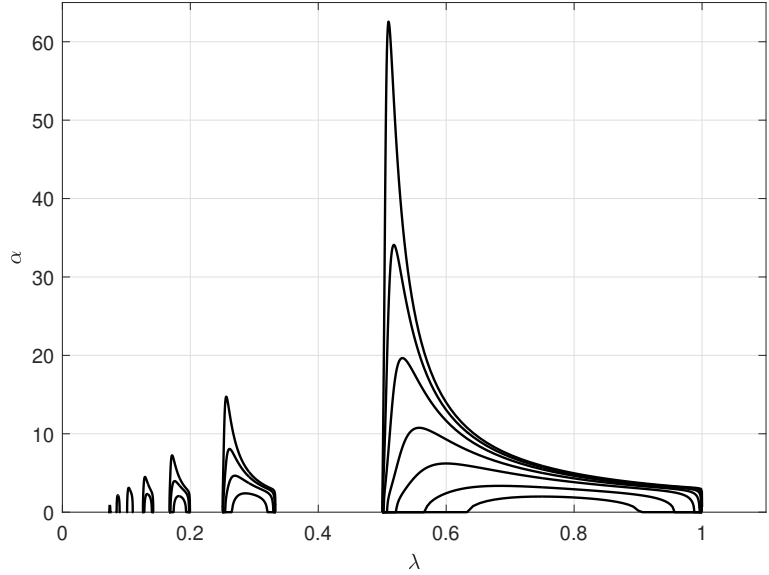}
\caption{The $(\lambda,\alpha)$ bifurcation diagram with $D = 3\times 10^{-6}$, $10^{-5}$, $3\times 10^{-5}$, $10^{-4}$, $3\times 10^{-4}$, $10^{-3}$, $2 \times 10^{-3}$. The amplitude, $\alpha$, increases as $D$ decreases.}
\label{fig:5.1}
\end{center}
\end{figure}

\begin{figure}
\begin{center}
\includegraphics[width=0.8\textwidth]{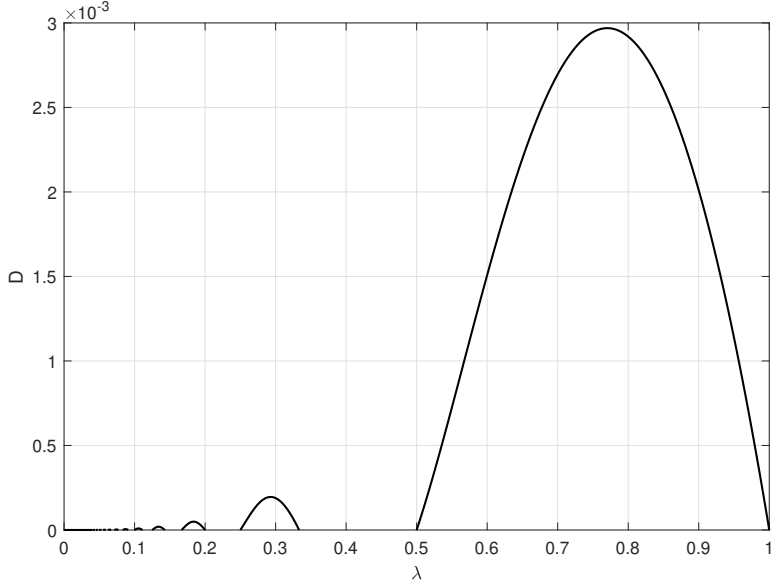}
\caption{The region $\Omega$ lies below these curves, or tongues.}
\label{fig:5.2}
\end{center}
\end{figure}

\begin{figure}
\begin{center}
\includegraphics[width=\textwidth]{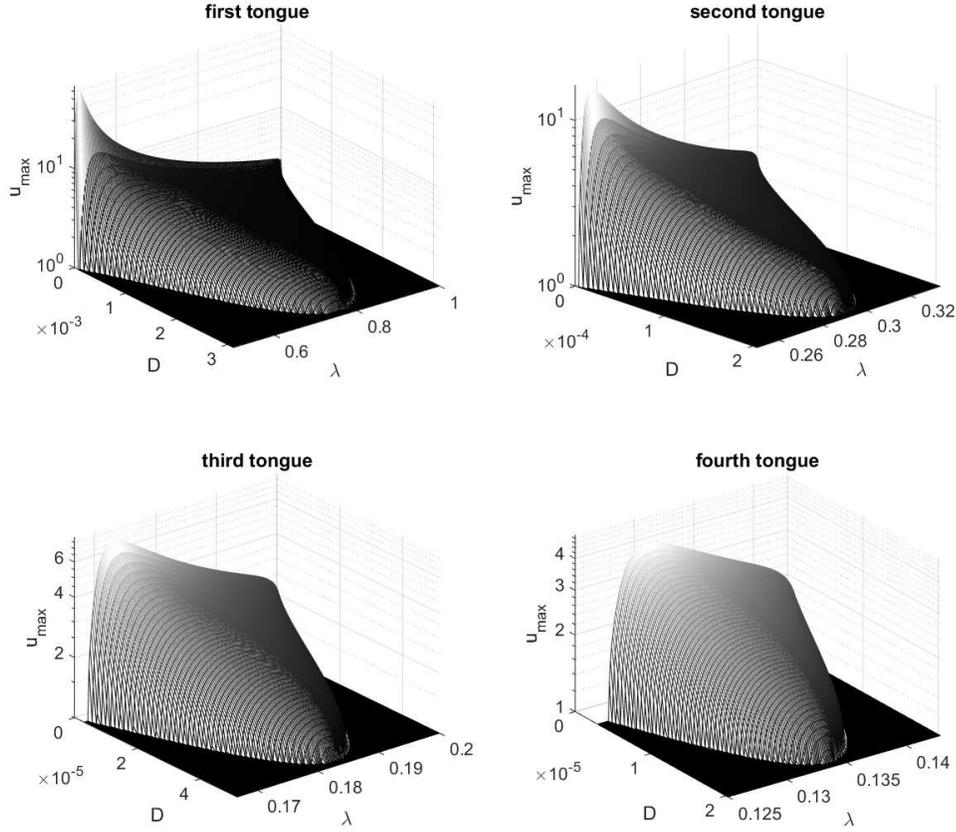}
\caption{Semilogarithmic plots of $u_{max}$ in the first four tongues of $\Omega$.}
\label{fig:5.3}
\end{center}
\end{figure}

The numerically-calculated $(\lambda,\alpha)$ bifurcation diagram is given in Figure \ref{fig:5.1} for various values of $D$.
It is now straightforward to construct the bifurcation locus in the $(\lambda,D)$ plane.
This is given directly by the intersection of the curve 
\begin{align}
    \notag D=-\frac{\lambda^3}{4\pi^3}\sin\frac{\pi}{\lambda}
\end{align}
with the positive quadrant of the $\lD$ plane.
This intersection is in the form of a countably infinite sequence of `tongue like' curves based on the $\lambda-$axis, and with monotone decreasing `height'.
The $i^{th}$ `tongue' has base points at $\lambda=1/(2i-1)$ and $\lambda=1/2i$, with `height' $D=\Delta_i$.
We label the set of interior points of the $i^{th}$ `tongue' as $\Omega_i$, and its boundary as $\partial \Omega_i=\overline{\Omega_i}\setminus\Omega_i $.
We write
\begin{align}
    \label{5.15} \Omega=\bigcup_{i=1}^\infty \Omega_i,\quad \partial\Omega=\bigcup_{i=1}^\infty \partial\Omega_i.
\end{align}
At each point $\lD\in\Omega$ there is a \emph{unique} (up to translation in $x$) and \emph{nontrivial} periodic steady state, with fundamental wavelength $\lambda$ and amplitude $\alpha(\lambda,D)$.
Since the equation \eqref{5.1} is invariant under the transformation $x\mapsto -x$, this allows for the existence of a translation in $x$, so that at each $\lD\in\Omega$ we may select a representative periodic steady state $u=F_p(x,\lambda,D)$ which is an \emph{even function} of $x\in\R$, so that 
\begin{align}
    \label{5.16} &F_p(-x,\lambda,D)=F_p(x,\lambda,D)\quad\forall x\in\R
    \\\label{5.17} &F_p'(0,\lambda,D)=0
\end{align}
and
\begin{align}
    \label{5.18} F_p'\left(\pm\frac{1}{2}\lambda,\lambda,D\right)=0.
\end{align}
A plot of $\Omega$ in the $\lD$ plane is shown in Figure~\ref{fig:5.2}.
The numerically-calculated maximum value of $u=F_p(x,\lambda,D)$, which we denote by $u_{max}$, is shown in Figure~\ref{fig:5.3} for the first four tongues of $\Omega$.  This is a more suitable measure of amplitude than $\alpha$ for the purposes of this semi-logarithmic plot, since $\alpha$ is zero at the edges of the region where the nonlinear solutions exist. 

We now examine the structure of the periodic steady states with $\lD\in\overline{\Omega}$.
Firstly, on the bifurcation locus, with $\lD\in\partial\Omega\cap(\R^+)^2$ we have 
\begin{align}
    \label{5.19} F_p(x,\lambda,D)=1\quad \forall x\in\R,
\end{align}
whilst it follows from \eqref{5.13} and \eqref{5.14} that 
\begin{align}
    \label{5.20} F_p(x,\lambda,D)=1+O(|\lambda-\lambda_0|^{\frac{1}{2}}+|D-D_0|^{\frac{1}{2}})\text{ as }(\lambda,D)\to(\lambda_0,D_0)
\end{align}
uniformly for $x\in\R$, with $(\lambda_0,D_0)\in \partial\Omega\cap(\R^+)^2$ and $\lD\in\Omega$. 
In terms of regularity, it follows from \eqref{5.2} and \eqref{5.1} (see for example, \cite{Hale1977}) that 
\begin{align}
    \label{5.21} F_p\in C^{\omega,0,0}(\R\times (\overline{\Omega}\cap (\R^+)^2)).
\end{align}
We can now establish the following bounds:
\begin{itemize}
    \item[(B1)] For any $\lD\in\Omega$, then, $$\inf_{x\in\R}F_p(x,\lambda,D)>0$$
\end{itemize}
\begin{proof}
Set $(\lambda^*,D^*)\in\Omega$, then $(\lambda^*,D^*)\in\Omega_i$ for some $i\in\N$. 
The point $(\lambda^*,D^*)$ can then be connected, in $\Omega_i$, to a point $(\lambda_0,D_0)\in\partial\Omega_i\cap (\R^+)^2$ via a straight line segment $L$.
Now, suppose that 
\begin{align}
    \label{5.22}\inf_{x\in\R}F_p(x,\lambda^*,D^*)\leq0.
\end{align}
We recall from \eqref{5.19} that
\begin{align}
    \label{5.23} \inf_{x\in\R}F_p(x,\lambda_0,D_0)=1.
\end{align}
Thus, from \eqref{5.21}, it follows that there exists $(\lambda_1,D_1)\in L\cap \Omega_i$, such that 
\begin{align}
    \label{5.24} \inf_{x\in\R} F_p(x,\lambda_1,D_1)=0.
\end{align}
As a consequence (since $F_p$ is periodic in $x\in\R$, together with \eqref{5.21}) there exists $x_1\in\R$ such that 
\begin{align}
    \label{5.25}F_p(x_1,\lambda_1,D_1)=0.
\end{align}
It then follows, via \eqref{5.24}, \eqref{5.25}, \eqref{5.21} and \eqref{5.1}, that 
\begin{align}
    \label{5.26} F_p^{(n)}(x_1,\lambda_1,D_1)=0\quad\forall\,n\in\N.
\end{align}
It is then a consequence of \eqref{5.26} with \eqref{5.21} that,
\begin{align}
    \notag F_p(x,\lambda_1,D_1)=0\quad\forall\,x\in\R.
\end{align}
However, $(\lambda_1,D_1)\in\Omega_i$, and we have established earlier in this section that as such $F_p(x,\lambda_1,D_1)$ must be a \emph{non-trivial} periodic function of $x\in\R$, with fundamental wavelength $\lambda_1$, and so cannot be the trivial function. Thus we arrive at a contradiction.
The result follows.
\end{proof}
\newcommand{\Fp}{{F_p(x,\lambda,D)}}
We next have
\begin{itemize}
    \item[(B2)] For any $\lD\in\Omega$, then, $$\sup_{x\in\R}F_p(x,\lambda,D)>1$$
\end{itemize}
\begin{proof}
Set $\lD\in\Omega$, then $\lD\in\Omega_i$ for some $i\in\N$.
Suppose that 
\begin{align}
    \label{5.27}\sup_{x\in\R}F_p(x,\lambda,D)\leq1.
\end{align}
Then since $F_p(x,\lambda,D)$ is a non-trivial periodic function in $x\in\R$, with fundamental wavelength $\lambda$, it follows that there exists $x^*\in\R$ such that 
\begin{align}
    \label{5.28} 0<F_p(x^*,\lambda,D)=\inf_{x\in\R}F_p(x,\lambda,D)<1, 
\end{align}
with the left hand inequality following from (B1).
From \eqref{5.27} and \eqref{5.28} it follows that,
\begin{align}
    \label{5.29}\int_{x^*-\frac{1}{2}}^{x^*+\frac{1}{2}}F_p(y,\lambda,D) dy<1.
\end{align}
Also, from \eqref{5.1}, we have 
\begin{align}
    \label{5.30}F_p''(x^*,\lambda,D)=-F_p(x^*,\lambda,D)\left(1-\int_{x^*-\frac{1}{2}}^{x^*+\frac{1}{2}}F_p(y,\lambda,D) dy\right)<0,
\end{align}
via \eqref{5.28} and \eqref{5.29}.
However, via \eqref{5.28} and \eqref{5.21}, $x=x^*$ is a local minimum point for $\Fp$, and so 
$$F_p''(x^*,\lambda,D)\geq 0$$
which contradicts \eqref{5.30}.
Thus,
$$\sup_{x\in\R}F_p(x,\lambda,D)>1,$$
as required.
\end{proof}
Correspondingly, we now have 
\begin{itemize}
    \item[(B3)] For any $\lD\in\Omega$, then, $$\inf_{x\in\R}F_p(x,\lambda,D)<1$$
\end{itemize}
\begin{proof}
This follows, with the obvious adjustments, that of (B2).
\end{proof}

To complete this section, we next consider in detail the structure of the periodic steady states in the $\lD$ plane as $D\to0^+$.
We begin by examining this limit for the periodic steady states in the principal tongue $\Omega_1$, and then give the corresponding results for $\Omega_i$, $i=2,3,4,\dots$ .

\subsection{Asymptotic structure as \texorpdfstring{$D\to0^+$}{text} with \texorpdfstring{$(\lambda,D)\in\Omega_1$}{text} }

We consider the asymptotic structure of $F_p(x,\lambda,D)$ with $(\lambda,D)\in\Omega_1$ as $D\to0^+$. This is achieved via a detailed application of the method of matched asymptotic expansions.
Due to the evenness of $F_p$ in $x$, we need only consider this structure for $x\in\left[0,\frac{1}{2}\lambda\right]$, with $\lambda\in\left(\frac{1}{2},1\right)$ fixed as $D\to0^+$.
We recall that, $F_p(x,\lambda,D)>0$ for all $x\in\left[0,\frac{1}{2}\lambda \right]$, whilst evenness requires 
\begin{align}
    \label{5.31} F_{p}^{\prime}(0, \lambda, D)=F_{p}^{\prime}\left(\frac{1}{2} \lambda, \lambda, D\right)=0.
\end{align}
An examination of \eqref{5.1}, together with a consideration of numerical solutions, dictates that, for $\lambda\in\left(\frac{1}{2},1\right)$ fixed, the asymptotic structure of $F_p(x,\lambda,D)$, as $D\to0^+$, develops into a three region structure as follows:

\noindent\emph{Region I} (support region) 
\begin{align}
    \notag x\in[0,S-O(D^\gamma)), \,\, F_p=O(1)^+\text{ as }D\to0^+
\end{align}
with $0<S<\frac{1}{2}\lambda$ and $\gamma>0$ to be determined. 
In general, $S$ will depend upon both $\lambda$ and $D$, with $S=O(1)^+$ as $D\to0^+$, and $\lambda\in\left(\frac{1}{2},1\right)$.
With this in mind, we expand $S(\lambda,D)$ as,
\newcommand{\as}{\text{ as }}
\begin{align}
    \label{5.32} S(\lambda,D)=\bar{S} (\lambda) +D^\gamma S_1(\lambda)+ o(D^\gamma),\as  D\to0^+
\end{align}
with $\lambda\in\left(\frac{1}{2},1\right)$, which allows for weak displacements in the location of region II below.

\noindent\emph{Region II} (transition region) 
\begin{align}
    \notag x\in\left(S-O(D^\gamma),S+O(D^\gamma)\right),\,\, F_p=O(D^\delta)^+\text{ as }D\to0^+
\end{align}
with $\delta>0$ to be determined. 

\noindent\emph{Region III} (exponential region)
\newcommand{\blank}{\underline{\hspace{1cm}}}
\begin{align}
    \notag x\in \left(S+O(D^\gamma),\frac{1}{2}\lambda\right],\,\, F_p=O(E(D))^+\as D\to0^+
\end{align}
with $E(D)$ indicating terms exponentially small in $D$ as $D\to0^+$.

The above structure requires, for the change in structure across region II, that the leading order interval length of region III (which is the leading order separation of consecutive support regions) must have the same length as the half span of the nonlocal term, which is $\frac{1}{2}$, and so,
\begin{align}
    \notag \bar{S}+\frac{1}{2}=\lambda-\bar{S}
\end{align}
which gives,
\begin{align}
    \label{5.33} \bar S =\frac{1}{2}\left(\lambda-\frac{1}{2}\right)=\bar S(\lambda)
\end{align}
and we note from this that 
\begin{align}
    \notag 0<\bar S(\lambda)<\frac{1}{4}
\end{align}
for $\lambda\in\left(\frac{1}{2},1\right)$, with,
\begin{align}
    \label{5.34} \bar{S}(\lambda)-\frac{1}{2}<-\bar{S}(\lambda)\quad  \forall\, \lambda \in\left(\frac{1}{2},1\right).
\end{align}
With \eqref{5.32}-\eqref{5.34}, then equation \eqref{5.1} requires that 
\begin{align}
    \label{5.35} \int_{x-\frac{1}{2}}^{x+\frac{1}{2}} F_p(y,\lambda,D)dy=\int_{-\frac{1}{2}}^{\frac{1}{2}} F_p(y,\lambda,D)dy+ O(E(D))\as D\to0^+,
\end{align}
for $x\in \left[0,\bar{S}(\lambda)-O(D^\gamma)\right)$, whilst 
\begin{align}
    \label{5.36} \int_{-\frac{1}{2}}^{\frac{1}{2}} F_p(y,\lambda,D)dy=1+\alpha_1 D+\beta_1 D^r +o(D^r)
\end{align}
as $D\to0^+$, with the constants $r>1$, $\alpha_1$ and $\beta_1$ to be determined.

We begin in region I and expand in the form
\begin{align}
    \label{5.37} F_p(x,\lambda,D)=F_0(x,\lambda)+D^m F_1(x,\lambda)+o(D^m)\as D\to0^+,
\end{align}
with $x\in[0,S(\lambda)-O(D^\gamma))$ as $m>0$ to be determined.
We substitute into \eqref{5.1}, to obtain, at leading order as $D\to0^+$,
\begin{align}
    \label{5.38} &F_0''-\alpha_1 F_0=0,\quad 0<x<\bar S(\lambda)
    \\\label{5.39} &F_0'(0,\lambda)=F_0(\bar S(\lambda),\lambda)=0
    \\\label{5.40} &F_0(x,\lambda)>0,\quad 0\leq x<\bar S(\lambda)
\end{align}
with the second condition in \eqref{5.39} required by asymptotic matching (via Van Dyke's principle, \cite{VanDykeMilton1975Pmif}) between region I and region II.
Finally, \eqref{5.36} requires
\begin{align}
    \label{5.41} \int_0^{\bar S(\lambda)}F_0(y,\lambda) dy=\frac{1}{2},
\end{align}
using the evenness of $F_p(x,\lambda,D)$, and its asymptotic form in regions II and III.
Now, the problem \eqref{5.38}-\eqref{5.40} is a classical Sturm-Liouville eigenvalue problem, with eigenvalue $-\alpha_1$, whilst condition \eqref{5.40} requires that $-\alpha_1$ is the \emph{lowest} eigenvalue.
This may be treated directly, to obtain, 
\begin{align}
    \label{5.42} \alpha_1=-\frac{\pi^2}{\left(\lambda-\frac{1}{2}\right)^2}
\end{align}
with
\begin{align}
    \label{5.43} F_0(x,\lambda)=B \cos\frac{\pi x}{\left(\lambda-\frac{1}{2}\right)},\quad 0\leq x\leq \bar S(\lambda)
\end{align}
and $B>0$ an arbitrary constant.
It remains to apply condition \eqref{5.41}, which gives,
\begin{align}
\label{5.44} B=\frac{\pi}{(2\lambda-1)},
\end{align}
and so,
\begin{align}
    \label{5.45} F_0(x,\lambda)=\frac{\pi}{(2\lambda-1)}\cos\frac{\pi x}{\left(\lambda-\frac{1}{2}\right)},\quad 0\leq x\leq \bar S(\lambda)
\end{align}
and $\bar{S}(\lambda)=\frac{1}{2}\left(\lambda-\frac{1}{2}\right)$. 
Via \eqref{5.35}, \eqref{5.36} and \eqref{5.42}, we also have
\begin{align}
    \label{5.46} \int_{x-\frac{1}{2}}^{x+\frac{1}{2}} F_p(y,\lambda,D)dy=2\int_0^{\frac{1}{2}}F_p(y,\lambda,D)dy=1-\frac{\pi^2}{\left(\lambda-\frac{1}{2}\right)^2}D+\beta_1 D^r+o(D^r)
\end{align}
as $D\to0^+$, with $x\in[0,S(\lambda)-O(D^\gamma))$.
This completes the leading order structure in region I.
We now move on to consider region II, the transition region.
In region II we introduce the scaled coordinate $X=O(1)$ as $D\to0^+$, given by,
\begin{align}
    \label{5.47} x=S(\lambda,D)+D^\gamma X
\end{align}
It then follows from the expansion in region I, on moving into region II, via \eqref{5.45}, \eqref{5.37} with \eqref{5.41}, that $F_p=O(D^\gamma)$ as $D\to0^+$ in region II.
Thus, we expand in the form
\begin{align}
    \label{5.48} F_p(X,\lambda,D)=D^\gamma \bar F_0(X,\lambda)+o(D^\gamma)\as D\to0^+
\end{align}
with $X=O(1)$. We remark here that the term at $O(D^\gamma)$ in $S(\lambda,D)$ could be removed by a constant shift in the coordinate $X$. However, it is convenient for later matching purposes to leave this shift in $S$; this also leads to a parameter free equation at leading order, with all parameters shifted into the matching condition.
We next consider the form of the nonlocal term in \eqref{5.1} in region II, which can be written as,
\begin{align}
    \notag &\int_{x-\frac{1}{2}}^{x+\frac{1}{2}} F_p(y,\lambda,D)dy =\int_{S(\lambda,D)-\frac{1}{2}+D^\gamma X}^{S(\lambda,D)+\frac{1}{2}+D^\gamma X} F_p(y,\lambda,D)dy
    \\\label{5.49} &=\int_{S(\lambda,D)-\frac{1}{2}+D^\gamma X}^{S(\lambda,D)-\frac{1}{2}} F_p(y,\lambda,D)dy+\int_{S(\lambda,D)-\frac{1}{2}}^{S(\lambda,D)+\frac{1}{2}-O(1)} F_p(y,\lambda,D)dy+\int_{S(\lambda,D)+\frac{1}{2}-O(1)}^{S(\lambda,D)+\frac{1}{2}+D^\gamma X} F_p(y,\lambda,D)dy.
\end{align}
Now, using the periodicity and evenness of $F_p$, and its structure in each of regions I-III, we have firstly that,
\begin{align}
\label{5.50} \int_{S(\lambda,D)-\frac{1}{2}+D^\gamma X}^{S(\lambda,D)-\frac{1}{2}} F_p(y,\lambda,D)dy=O(E(D)D^\gamma)
\end{align}
as $D\to0^+$ with $X=O(1)$, whilst,
\begin{align}
    \label{5.51} \int_{S(\lambda,D)-\frac{1}{2}}^{S(\lambda,D)+\frac{1}{2}-O(1)} F_p(y,\lambda,D)dy=1-\frac{\pi^2}{\left(\lambda-\frac{1}{2}\right)^2}D+\beta_1 D^r+o(D^r)
\end{align}
as $D\to0^+$.
The final term in \eqref{5.49} may now be written as
\begin{align}
    \notag &\int_{S(\lambda,D)+\frac{1}{2}-O(1)}^{S(\lambda,D)+\frac{1}{2}+D^\gamma X} F_p(y,\lambda,D)dy = D^\gamma\int_{Y=-O(D^{-\gamma})}^{Y=X} F_p\left(S(\lambda,D)+\frac{1}{2}+D^\gamma Y,\lambda,D\right)dY
    \\\notag &= D^\gamma\int_{Y=-O(D^{-\gamma})}^{Y=X} F_p(\lambda-S(\lambda,D)+D^\gamma Y,\lambda,D)dY 
    \\&= D^\gamma\int_{Y=-O(D^{-\gamma})}^{Y=X} F_p(-S(\lambda,D)+D^\gamma Y,\lambda,D)dY\label{5.52}
    \\ \notag&= D^\gamma\int_{Y=-O(D^{-\gamma})}^{Y=X} F_p(S(\lambda,D)-D^\gamma Y,\lambda,D)dY \sim D^{2\gamma}\int_{Y=-\infty}^{Y=X} \bar F_0(-Y,\lambda)dY \\\notag&=D^{2\gamma}\int_{z=-X}^{z=\infty} \bar F_0(z,\lambda)dz
\end{align}
as $D\to0^+$ with $X=O(1)$.
Using \eqref{5.50}-\eqref{5.52} in \eqref{5.49}, we finally arrive at
\[
    \int_{S(\lambda,D)-\frac{1}{2}+D^\gamma X}^{S(\lambda,D)+\frac{1}{2}+D^\gamma X} F_p(y,\lambda,D)dy\]
\begin{equation}
    =1-\frac{\pi^2}{\left(\lambda-\frac{1}{2}\right)^2}D+\beta_1 D^r+D^{2\gamma}\int_{z=-X}^{z=\infty} \bar F_0(z,\lambda)dz+o(D^{2\gamma},D^r,E(D))\label{5.53} 
\end{equation}
as $D\to0^+$ with $X=O(1)$.
We now substitute from \eqref{5.47}, \eqref{5.48} and \eqref{5.53} into equation \eqref{5.1}.
To obtain a non-trivial balance at leading order requires us to choose 
\begin{align}
    \label{5.54} \gamma=\frac{1}{4}.
\end{align}
At leading order, equation \eqref{5.1} then becomes,
\begin{align}
    \label{5.55} \Psi_{\bar X \bar X}-\Psi \int_{w=-\bar X}^{w=\infty}\Psi(w)dw=0,\quad \bar X\in\R,
\end{align}
which is both nonlinear and nonlocal.
In \eqref{5.55} we have introduced the simple scalings
\begin{align}
    \label{5.56}
        \bar F_0=(2\lambda-1)^{-\frac{3}{2}}\Psi,~~~X =(2\lambda -1)^{\frac{1}{2}}\bar X.
\end{align}
Equation \eqref{5.55} is now completed by matching conditions between region II and region I (as $\bar{X}\to -\infty$) and region III (as $\bar{X}\to\infty$).
The matching process is straightforward (following Van Dyke's asymptotic matching principle, \cite{VanDykeMilton1975Pmif}), and leads to the two boundary conditions,
\begin{align}
    \label{5.57} &\Psi(\bar X)\to0\as \bar X\to\infty,
    \\\label{5.58} &\Psi(\bar X)\sim -2\pi^2(\bar X+l)\as \bar X\to-\infty.
\end{align}
Here
\begin{align}
    \label{5.59} S_1(\lambda)=(2\lambda-1)^{\frac{1}{2}} l.
\end{align}
It is also required that,
\begin{align}
    \label{5.60} \Psi(\bar X)>0\quad\forall\,\bar X\in \R.
\end{align}

\begin{figure}
\begin{center}
\includegraphics[width=0.8\textwidth]{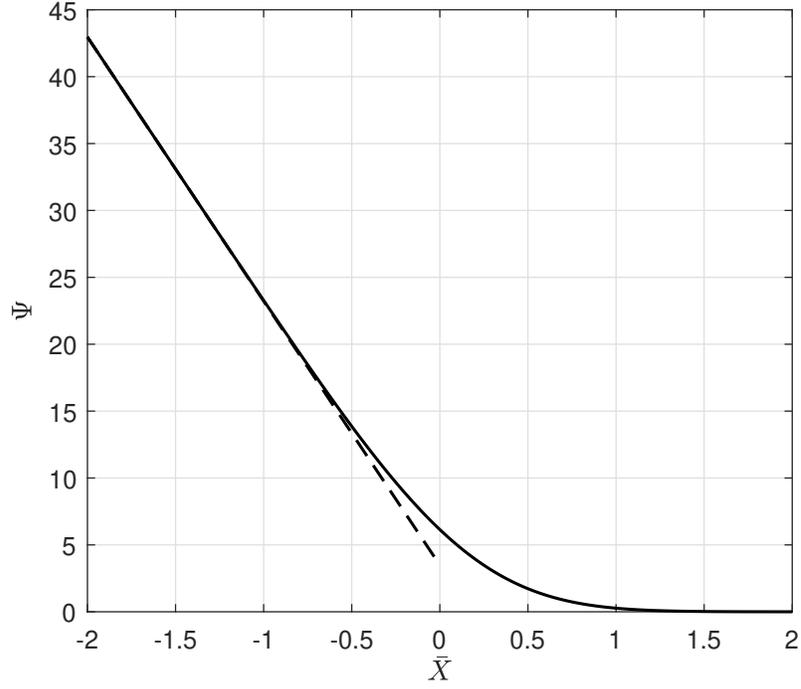}
\caption{The numerically-calculated solution of (\ref{5.55}) subject to (\ref{5.56}). The broken line is $\Psi = -2 \pi^2 (\bar{X}+l^*)$.}
\label{fig:5.4}
\end{center}
\end{figure}

We observe that the problem \eqref{5.55}, \eqref{5.57}, \eqref{5.58} and \eqref{5.60} is a nonlinear, nonlocal eigenvalue problem, with real eigenvalue $l$. It is also free of parameters.
A numerical investigation of this problem establishes that a solution exists if and only if $l=l^*$, and the solution is unique.
Here
\begin{align}
    \label{5.61} l^*\approx-\frac{3.493}{2\pi^2},
\end{align}
and a numerically-determined graph of $\Psi(\bar{X})$ against $\bar{X}$, calculated using a simple finite difference method and the trapezium rule, is shown in Figure~\ref{fig:5.4}.
It is also straightforward to establish that
\begin{align}
    \label{5.62} \Psi(\bar X)=-2\pi^2(\bar X+l^*)+\Psi_{-\infty}\bar X^{-2}e^{-\frac{1}{2}\pi X^2}(1+o(1))
\end{align}
as $\bar{X}\to -\infty$, whilst,
\begin{align}
    \label{5.63} \Psi(\bar X)=\Psi_{\infty} e^{-\frac{1}{2}\pi X^2}(1+o(1))
\end{align}
as $\bar{X} \to+\infty$. 
Here $\Psi_\infty$ ($>0$) and $\Psi_{-\infty}$ are parameter free, globally determined constants. In principle, these can be determined from the numerical solution shown in Figure~\ref{fig:5.4}, but in practice, since these corrections are exponentially small, we are not able to resolve them with sufficient accuracy.

We now return to region I. To allow matching with region II we now require 
\begin{align}
    \label{5.65} m=\frac{1}{4},
\end{align}
whilst the most structured balance in equation \eqref{5.1} dictates that 
\begin{align}
    \label{5.66}r=m+1=\frac{5}{4}.
\end{align}
The problem for $F_1$ is then,
\begin{align}
    \label{5.67} F_{1}''+\frac{\pi^{2}}{\left(\lambda-\frac{1}{2}\right)^{2}} F_{1}&=\beta_1 F_{0}(x, \lambda)
    \\\notag &=\frac{\beta_1\pi}{(2\lambda-1)}\cos\frac{\pi x}{\left(\lambda-\frac{1}{2}\right)},\quad 0<x<\frac{1}{2}\left(\lambda -\frac{1}{2}\right), 
\end{align}
\begin{align}
    \label{5.68}&F_1' (0,\lambda)=F_1\left(\frac{1}{2}\left(\lambda-\frac{1}{2}\right),0\right)=0,
    \\\label{5.69}&\int_0^{\frac{1}{2}\left(\lambda-\frac{1}{2}\right)} F_1(y,\lambda)dy=0.
\end{align}
The second of conditions \eqref{5.68} arises from matching with region II (with both expansions taken up to $O(D^{\frac{1}{4}})$ in regions I and II).
A solution to \eqref{5.67} and \eqref{5.68} exists if and only if 
\begin{align}
    \label{5.70} \beta_1=0,
\end{align}
after which
\begin{align}
    \label{5.71}F_1(x,\lambda)=C\cos\frac{\pi x}{\left(\lambda-\frac{1}{2}\right)}
\end{align}
with $C$ an arbitrary real constant. 
However, substitution from \eqref{5.71} into \eqref{5.69}, then requires $C=0$, and so,
\begin{align}
    \label{5.72} F_1(x,\lambda)=0\quad\forall\,x\in\left[0,\frac{1}{2}\left(\lambda-\frac{1}{2}\right)\right].
\end{align}
Thus, in region I we have 
\begin{align}
    \label{5.73} F_p(x,\lambda,D)=\frac{\pi}{(2\lambda-1)}\cos\frac{\pi x}{\left(\lambda-\frac{1}{2}\right)}+o(D^{\frac{1}{4}})
\end{align}
as $D\to0^+$ with $x\in\left[0,\frac{1}{2}\left(\lambda-\frac{1}{2}\right)\right)$.

We finally move into the exponential region, region III.
In this region $F_p$ is exponentially small in $D$ as $D\to0^+$, and develops as a WKB expansion.
For brevity we omit details, but obtain
\begin{align}
    \label{5.74} F_p(x,\lambda,D)=A^*(\lambda)e^{-\Phi_0(\lambda)D^{-\frac{1}{2}}}\cosh(\frac{\Phi(x)}{D^\frac{1}{2}})(1+o(1))
\end{align}
as $D\to0^+$, with $x\in\left(\frac{1}{2}(\lambda -\frac{1}{2})+O(D^{\frac{1}{4}}),\frac{1}{2}\lambda\right]$.
Here,
\begin{align}
    \label{5.75}\Phi(x)=\frac{1}{\sqrt 2}\int_x^{\frac{1}{2}\lambda}\left(1-\sin(\frac{\pi w}{\left(\lambda-\frac{1}{2}\right)})\right)^{\frac{1}{2}}dw\quad \forall\,x\in\left[\frac{1}{2}\left(\lambda-\frac{1}{2}\right),\frac{1}{2}\lambda\right)
\end{align}
and
\begin{align}
    \label{5.76}\Phi_0(\lambda)=\Phi\left(\frac{1}{2}\left(\lambda-\frac{1}{2}\right)\right).
\end{align}
The constant $A^*(\lambda)$ is independent of $D$ and is related to the constant $\Psi_\infty$ and $\lambda$, via matching expansion \eqref{5.74} as $x\to\frac{1}{2}\left(\lambda-\frac{1}{2}\right)^-$ with expansion \eqref{5.48} as $X\to\infty$, which gives,
\begin{align}
    \label{5.77}A^*(\lambda)=\frac{2\Psi_\infty}{(2\lambda-1)^{\frac{3}{2}}}
\end{align}
with details omitted for brevity. This completes the asymptotic structure as $D\to0^+$, with fixed $\lambda\in\left(\frac{1}{2},1\right)$. 

\begin{figure}
\begin{center}
\includegraphics[width=0.8\textwidth]{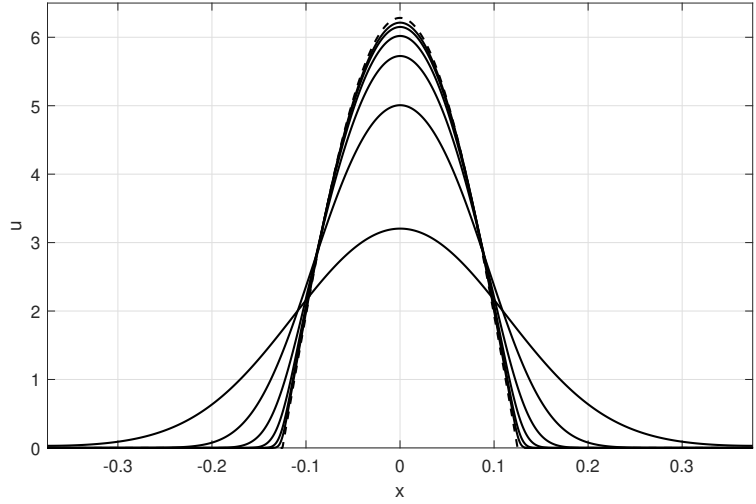}
\caption{The solution $F_p$, calculated numerically for $\lambda = \frac{3}{4}$ and $D = 10^{-3}$, $10^{-4}$, $10^{-5}$, $10^{-6}$, $10^{-7}$ and $10^{-8}$. The broken line shows the leading order outer solution given by (\ref{5.73}).}
\label{fig:5.5}
\end{center}
\end{figure}

We make the following key observations. First, up to terms exponentially small in D, as $D\to0^+$,
\begin{align}
    \label{5.78} \supp(F_p(x,\lambda,D))=[-S(\lambda,D),S(\lambda,D)]
\end{align}
with
\begin{align}
    \label{5.79} S(\lambda,D)=\frac{1}{2}\left(\lambda-\frac{1}{2}\right)-\frac{3.493}{\sqrt 2\pi^2}\left(\lambda-\frac{1}{2}\right)^{\frac{1}{2}}D^{\frac{1}{4}}+o(D^{\frac{1}{4}})\,\as D\to0^+, 
\end{align}
whilst,
\begin{align}
    \label{5.80}\alpha (\lambda,D)=\frac{\pi}{(2\lambda-1)}+o(D^{\frac{1}{4}})\,\as D\to0^+
\end{align}
and occurs at $x=0$. In addition, 
\begin{align}
    \label{5.81}\int_{-S(\lambda,D)}^{+S(\lambda,D)} F_p(y,\lambda,D)dy=1-\frac{\pi^2}{\left(\lambda-\frac{1}{2}\right)^2}D+o(D^{\frac{5}{4}})
\end{align}
as $D\to0^+$. 
For comparison, we present in Figure~\ref{fig:5.5}, a numerically determined representation of $F_p$, with $\lambda=0.75$ and various values of $D$.
This is in excellent agreement with the asymptotic form developed above.

Finally, we observe from \eqref{5.78}-\eqref{5.81} in particular, that the asymptotic structure developed above, for $F_p(x,\lambda,D)$ as $D\to0^+$, with fixed $\lambda\in \left(\frac{1}{2},1\right)$, becomes nonuniform as $\lambda\to\frac{1}{2}$, and, specifically, when $\lambda=\frac{1}{2}+O(D^{\frac{1}{2}})$ as $D\to0^+$.
For the present, we continue by considering the asymptotic structure of $F_p(x,\lambda,D)$ as $D\to0^+$ in each of the remaining tongues $\Omega_i$, $i=2,3,\dots$ .

\subsection{Asymptotic structure as \texorpdfstring{$D\to0^+$}{text} with \texorpdfstring{$(\lambda,D)\in\Omega_i$}{text}, \texorpdfstring{$i=2,3,\dots$}{text}}

We now consider the asymptotic form of $F_p(x,\lambda,D)$ with $\lD\in\Omega_i$ ($i\geq2$) as $D\to0^+$.
We, again, need only consider the structure for $x\in\left[0,\frac{1}{2}\lambda\right]$, with now $\lambda\in\left(\frac{1}{2i},\frac{1}{2i-1}\right)$ fixed, as $D\to0^+$.
The structure is again as in the case of $\Omega_1$, and for brevity, we restrict attention to the key support region, labelled as region I.
The analysis follows exactly that of the previous subsection, and we therefore present only the salient features.
In this case we obtain,
\begin{align}
    \label{5.82}&\bar S (\lambda)= \frac{1}{2}i\left(\lambda-\frac{1}{2i}\right),
    \\\label{5.83}& F_0(x,\lambda)=\frac{\pi}{(2i-1)(2i\lambda-1)}\cos\frac{\pi x}{i\left(\lambda-\frac{1}{2i}\right)}\quad\forall\,x\in[0,\bar S(\lambda)],
    \\\label{5.84}& (2i-1)\int_{-S(\lambda,D)}^{S(\lambda,D)}F_p(y,\lambda,D)dy=1-\frac{\pi^2}{i^2\left(\lambda-\frac{1}{2i}\right)^2}D+o(D)
\end{align}
as $D\to0^+$, whilst,
\begin{align}
    \label{5.85} \alpha(\lambda,D)=\frac{\pi}{(2i-1)(2i\lambda-1)}+o(1)\as D\to0^+
\end{align}
with $\lambda \in\left(\frac{1}{2i},\frac{1}{2i-1}\right)$ fixed.
We observe that, as before, this structure becomes nonuniform when $\lambda=\frac{1}{2i}+o(1)$ as $D\to0^+$.

\subsection{Asymptotic structure as \texorpdfstring{$D\to0^+$}{text} with \texorpdfstring{$\lambda=\frac{1}{2} +O(D^{\frac{1}{2}})$}{text}}\label{sec_5.3.3}
\newcommand{\lambdabar}{\bar{\lambda}}
\newcommand{\xhat}{\hat{X}}
We return to $(\lambda, D)\in\Omega_1$, and now develop the structure to $F_p(x,\lambda,D)$ with $\lambda=\frac{1}{2} +O(D^{\frac{1}{2}})$ and $x\in\left[0,\frac{1}{2}\lambda\right]$ as $D\to0^+$.
We write
\begin{align}
    \label{5.86}\lambda=\frac{1}{2}+D^{\frac{1}{2}}\bar \lambda
\end{align}
with $\bar{\lambda}=O(1)^+$ as $D\to0^+$.
An examination of \eqref{5.78}-\eqref{5.80}, with \eqref{5.86}, then requires,
\begin{align}
    \label{5.87}x=O(D^{\frac{1}{2}})^+ \text{ and } F_p=O(D^{-\frac{1}{2}})^+
\end{align}
as $D\to0^+$ in the support region, with 
\begin{align}
   \label{5.88} F_p=O(E(D)),
\end{align}
when 
\begin{align}
    \label{5.89} x\in\left(O(D^{\frac{1}{2}}),\frac{1}{4}+\frac{1}{2}\bar{\lambda} D^{\frac{1}{2}}\right]
\end{align}
in the exponential region.
Following \eqref{5.87}, in the support region we introduce
\begin{align}
    \label{5.90} \xhat=\frac{x}{D^{\frac{1}{2}}}=O(1)^+
\end{align}
as $D\to0^+$, and expand in the form
\begin{align}
    \label{5.91}F_p(\xhat,\lambdabar,D)=v(\xhat,\lambdabar)D^{-\frac{1}{2}}+o(D^{-\frac{1}{2}})
\end{align}
as $D\to0^+$ with $\xhat=O(1)^+$.
After a careful consideration of the nonlocal term, substitution from \eqref{5.86}, \eqref{5.90} and \eqref{5.91} into \eqref{5.1} gives, at leading order, the nonlocal, nonlinear equation,
\begin{align}
    \label{5.92} v_{\xhat \xhat} +v\left(1-2\int_{-\infty}^\infty v(s,\lambdabar)ds +\int_{\xhat-\lambdabar}^{\xhat+\lambdabar}v(s,\lambdabar)ds \right)=0;\quad \xhat>0
\end{align}
which must be solved subject to the conditions,
\begin{align}
    \label{5.93}&v_{\xhat}(0,\lambdabar)=0,
    \\\label{5.94}&v(\xhat,\lambdabar)\to0\as \xhat\to\infty,
    \\\label{5.95}&v(\xhat,\lambdabar)>0 \quad \forall\, x\geq0,
    \\\label{5.96}& v(-\xhat,\lambdabar)=v(\xhat,\lambdabar) \quad \forall\xhat\geq0.
\end{align}
Here, the boundary condition \eqref{5.94} accommodates asymptotic matching to the exponential region. 

\begin{figure}
\begin{center}
\includegraphics[width=0.8\textwidth]{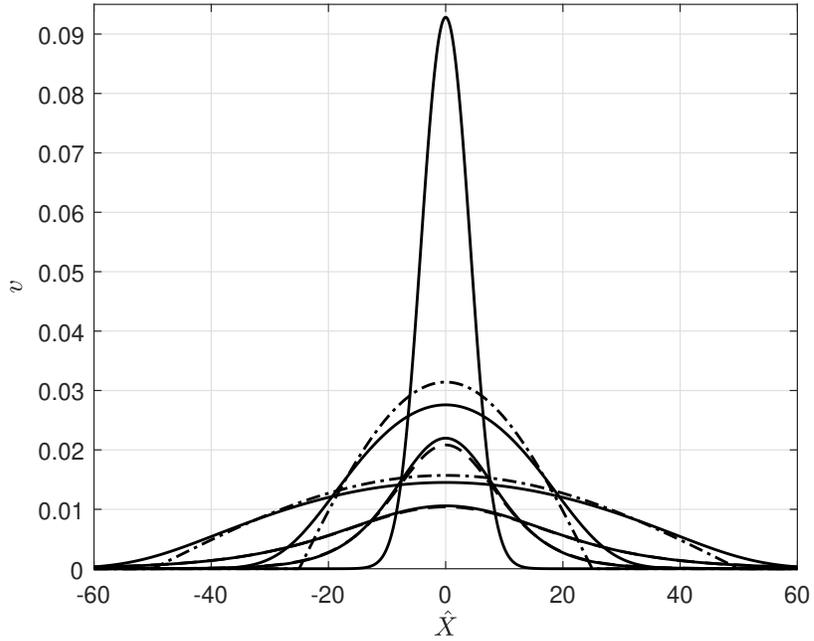}
\caption{The numerical solution (bold curves) of (\ref{5.92}) to (\ref{5.96}) for $\bar{\lambda} = 0.5$, $1$, $10$, $50$ and $100$. The broken line is the asymptotic solution for $\bar{\lambda}\ll 1$ for $\lambdabar = 0.5$ and $1$, given by (\ref{5.111}), whilst the dash-dotted line is the asymptotic solution for $\lambdabar \gg 1$ for $\lambdabar = 50$ and $100$, which comes from rescaling (\ref{5.73}).}
\label{fig:5.6}
\end{center}
\end{figure}

\begin{figure}
\begin{center}
\includegraphics[width=\textwidth]{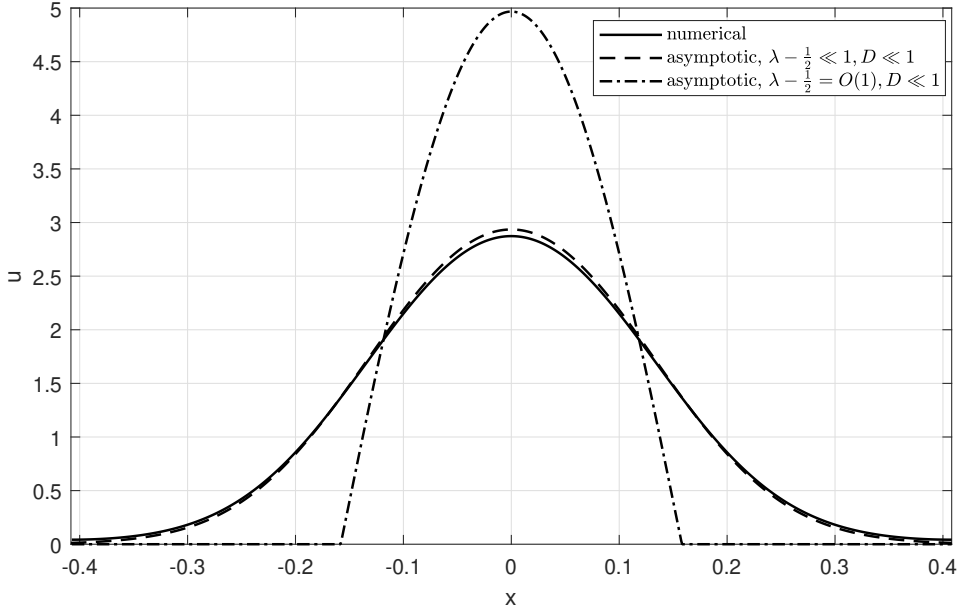}
\caption{The numerically calculated periodic solution of the full problem for $D = 10^{-3}$ and $\lambda = \frac{1}{2} + 10 \sqrt{D}$, along with the corresponding asymptotic solution for $\bar{\lambda} = 10$ (broken line) and the solution (\ref{5.73}), (dash-dotted line).}
\label{fig:5.7}
\end{center}
\end{figure}

\begin{figure}
\begin{center}
\includegraphics[width=\textwidth]{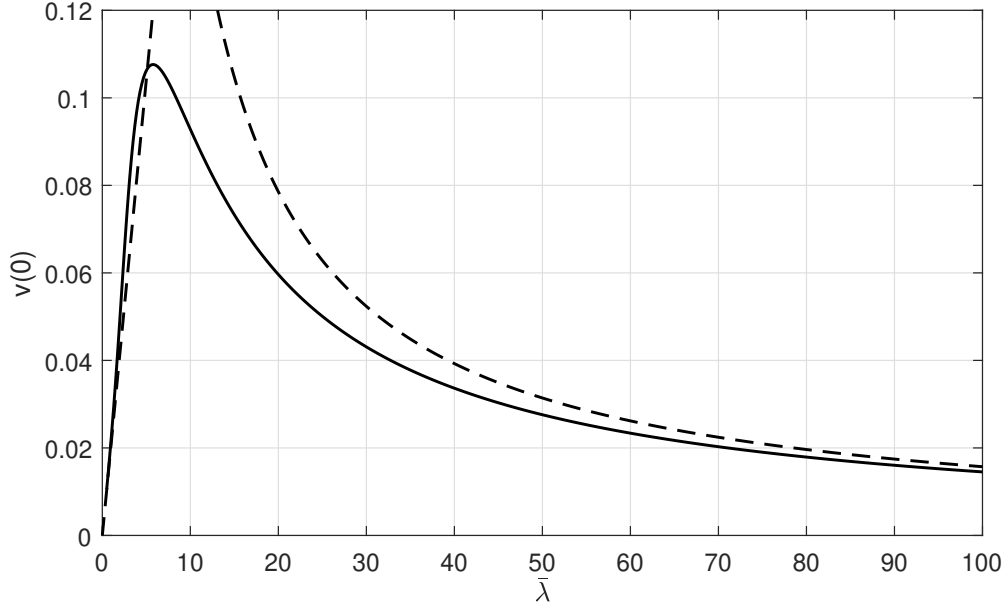}
\caption{The maximum value, $v(0,\bar{\lambda})$ of the numerically-calculated solution for  (\ref{5.92}) to (\ref{5.96}). The broken lines show the predicted asymptotic behaviour as $\bar{\lambda} \to 0$ and $\bar{\lambda} \to \infty$.}
\label{fig:5.8}
\end{center}
\end{figure}

The problem \eqref{5.92}-\eqref{5.96} has been considered numerically, which has demonstrated that a unique solution exists for each $\lambdabar>0$.
We observe that this solution has,
\begin{align}
    \label{5.97} v(\xhat,\lambdabar)\sim v_\infty (\lambdabar) e^{-\sigma_\infty(\lambdabar)\xhat}\as \xhat \to\infty
\end{align}
with
\begin{align}
    \label{5.98} \sigma_\infty(\lambdabar)=4\int_0^\infty v(s,\lambdabar)ds-1>0
\end{align}
and $v_\infty(\lambdabar)>0$ a globally dependent constant. 
Figure~\ref{fig:5.6} shows the numerical solution of (\ref{5.92}) to (\ref{5.96}) for various values of $\lambdabar$. The solution becomes larger and narrower as $\lambdabar$ increases, but then becomes smaller and fatter as $\lambdabar$ increases past about five (see also Figure~\ref{fig:5.8}).  Also shown is a comparison with the asymptotic solutions for $\lambdabar$ both large and small, with which there is excellent agreement.
Figure~\ref{fig:5.7} shows the numerically calculated periodic solution of the full problem for $D = 10^{-3}$ and $\lambda = \frac{1}{2} + 10 \sqrt{D}$, along with the corresponding asymptotic solution for $\bar{\lambda} = 10$ and the solution (\ref{5.73}).
The agreement is excellent, and we can also see that the $\lambda = O(1)$ asymptotic solution, (\ref{5.73}), is not a good approximation to the solution if $\lambda - \frac{1}{2}$ is of $O(\sqrt{D})$, as expected. The maximum of $v$ is at $\xhat=0$ and a graph of $v(0,\lambdabar)$ against $\lambdabar$ is shown in Figure~\ref{fig:5.8}, which shows that there is a single maximum at $\lambdabar \approx 5.8$, with, as required,
\begin{align}
    \label{5.99}\begin{split}
        v(0,\lambdabar)&\sim \frac{\pi}{2\lambdabar}\as \lambdabar\to\infty,
        \\v(0,\lambdabar)&\to0^+\as \lambdabar\to 0^+
    \end{split}
\end{align}
We note that,
\begin{align}
    \label{5.100}\alpha (\lambdabar,D)=D^{-\frac{1}{2}}v(0,\lambdabar)+o(D^{-\frac{1}{2}})
\end{align}
as $D\to0^+$. The amplitude given by $\alpha(\lambdabar,D)$ is indistinguishable from the amplitude given by the numerical solution of the full problem, shown in Figure~\ref{fig:5.1}, for $D$ less than about $10^{-4}$.

To complete this subsection, we examine the problem \eqref{5.92}-\eqref{5.96} in the limit $\lambdabar\to0^+$.
A balancing of terms in \eqref{5.92} leads us to write,
\newcommand{\xtilde}{\Tilde{X}}
\newcommand{\Itilde}{\Tilde{I}}
\newcommand{\vtilde}{\Tilde{v}}
\begin{align}
    \label{5.101} \xtilde=\lambdabar\xhat=O(1)^+\as \lambdabar\to0^+,
\end{align}
and expand in the form
\begin{align}
    \label{5.102} v(\xtilde,\lambdabar)=\lambdabar \vtilde(\xtilde)+o(\lambdabar) \as \lambdabar\to0^+.
\end{align}
Using \eqref{5.101} and \eqref{5.102}, the nonlocal term in \eqref{5.92} becomes 
\begin{align}
    \int_{\xhat-\lambdabar}^{\xhat+\lambdabar}v(s,\lambdabar)ds =\int_{\xhat-\lambdabar}^{\xhat+\lambdabar}(\lambdabar \vtilde(\lambdabar s)+o(\lambdabar))ds =\int_{\lambdabar\xhat-\lambdabar^2}^{\lambdabar\xhat+\lambdabar^2}(\vtilde(w)+o(1))dw =2\lambdabar^2\vtilde(\xtilde)+o(\lambdabar^2)\label{5.103}
\end{align}
as $\lambdabar\to0^+$. The most structured balance at leading order then requires,
\begin{align}
    \label{5.104}\int_{-\infty}^\infty v(\xtilde,\lambdabar)d\xtilde=\frac{1}{2}\lambdabar+\Itilde\lambdabar^3+o(\lambdabar^3)
\end{align}
as $\lambdabar\to0^+$ with $\Itilde$ a constant to be determined.
On substitution from \eqref{5.101}-\eqref{5.104} into \eqref{5.92}-\eqref{5.96}, we arrive at the leading order problem,
\begin{align}
    \label{5.105}&\vtilde_{\xtilde\xtilde}-2\vtilde(\Itilde-\vtilde)=0,\, \xtilde>0,
    \\\label{5.106}&\vtilde_{\xtilde}(0)=0,
    \\\label{5.107}&\vtilde(\xtilde)\to 0\as \xtilde\to\infty,
    \\\label{5.108}&\vtilde(\xtilde)>0\quad \xtilde\geq 0,
    \\\label{5.109}&\int_0^\infty \vtilde(w)dw=\frac{1}{4}.
\end{align}
This problem is a nonlinear eigenvalue problem, with eigenvalue $\Itilde\in\R$, which can be solved directly.
There is a single eigenvalue given by
\begin{align}
    \label{5.110}\Itilde=\frac{1}{72}
\end{align}
with the associated eigenfunction uniquely determined as
\begin{align}
    \label{5.111}\vtilde(\xtilde)=\frac{1}{48}\sech^2\left(\frac{1}{12}\xtilde\right)\quad \forall~\xtilde\geq 0.
\end{align}
Thus, \eqref{5.102} and \eqref{5.111} give
\begin{align}
    \label{5.112}v(0,\lambdabar)\sim\frac{1}{48}\lambdabar\,\as \lambdabar\to0^+
\end{align}
which is shown in Figure~\ref{fig:5.8}, and \eqref{5.109} and \eqref{5.110} lead to 
\begin{align}
    \label{5.113}\int_{-\infty}^\infty v(\xhat,\lambdabar)d\xhat=\frac{1}{2}+\frac{1}{72}\lambdabar^2 +o(\lambdabar^2)
\end{align}
as $\lambdabar\to0^+$.
Equations \eqref{5.111}-\eqref{5.113} are in excellent agreement with the numerical solution to problem \eqref{5.92}-\eqref{5.96}, when $\lambdabar$ is small

A final remark in this section relates to the temporal stability of the periodic steady states identified with $\lD\in\Omega$.
In this respect, we have performed a limited stability analysis of $u=F_p(x,\lambda,D)$ with $\lD\in\Omega$ when $D$ is small by considering the exponentially small part of the solution. Without presenting details, this analysis reveals that each such periodic steady state is locally, temporally, asymptotically stable to perturbations in the region between the $O(1)$ parts of the solution, where we would expect any instability to occur (extended intervals where $u$ remains uniformly very small, when long enough, can destabalise as they inherit the linear instability characteristics of the zero equilibrium). In relation to (IBVP), we observe that the results of this section support the possibility (P2), for the case when $0<D<\Delta_1$.\\
We readily verify that the wavelength selected by (IBVP) has, \emph{when $D$ is close to $\Delta_1$},
\begin{align}
    \label{5.114}\lambda\approx\lambda_m(D)=\frac{2\pi}{k_m(D)}\in \Omega_1
\end{align}
 and so we have, in (P2),
\begin{align}
    P(x) = F_p(x,\lambda_m(D),D),~ x\in \mathbb{R},
\end{align}
However this wavelength decreases as $D$ decreases, \emph{approaching the limiting value of $\frac{1}{2}$ as $D\to 0^+$}, and so now in (P2),
\begin{align}
    P(x) = F_p(x,\lambda(D),D),~ x\in \mathbb{R},
\end{align}
\emph{where now $\lambda(D)\to \frac{1}{2}^+$ as $D\to 0^+$}, which therefore has the structure given in detail in section 4.3. In particular the periodic solution selected therefore is of the nature of a spike, with base width of $O(D^{\frac{1}{2}})$ and height of $O(D^{-\frac{1}{2}})$. This is in full accord with the numerical solution of section 2.
A final point to note is that the detailed analysis of the existence and structure of the family of steady periodic solution in the first tongue $\Omega_1$ on the $(\lambda,D)$ plane given in this section not only provides substance to the conjecture (P2) in the present context, but also  plays a crucial role in studying the corresponding Cauchy problem on a closed finite spatial interval (with either Dirichlet or Neumann endpoint conditions), which we present in the second of this series of papers, \cite{BN2023_1}. 

We now move on to consider bifurcations to periodic travelling waves from the equilibrium solution $u_e=1$. This enables us to investigate in more depth the structure of the solution to (IBVP) when $t$ is large, now in the region of the propagating wavefronts identified in section 4. In particular, the next section first rules out the possibility of \emph{nondecaying, permanent form periodic waves travelling with, and to the rear of, the principal wavefronts.}

\section{Positive Periodic Travelling Waves}\label{section6}

We seek to identify the bifurcation to \emph{spatially periodic travelling wave} solutions from the equilibrium state $u_e=1$.
At fixed $D>0$, a periodic travelling wave solution, with propagation speed $v=v_b\neq 0$ will bifurcate from $u_e=1$, at wavelength $\lambda_b>0$  if and only if, $(\lambda_b,v_b)$ satisfies the complex equation,
\begin{align}
    \label{6.1}4 \pi^{2} D-2 i \pi v_{b} \lambda_{b}+\frac{\lambda_{b}^{3}}{\pi} \sin \frac{\pi}{\lambda_{b}} = 0
\end{align}
which is readily obtained from equation \eqref{1.3} with \eqref{1.9}, when written in the travelling wave coordinate $z=x-v t$, and linearized about $u_e=1$.
It is immediate that \eqref{6.1} has \emph{no solutions} which have $v_b\neq 0$ and $\lambda_b>0$.
We can therefore conclude that there are no bifurcations to periodic travelling waves from the equilibrium state $u_e=1$.\\
We are now able to use this result in moving on to consider transitional permanent form travelling waves, which accommodate the transition from the equilibrium state $u=0$ ahead of the wavefront to the equilibrium state $u=1$ or to a spatially periodic permanent form travelling wave, to the rear of the wavefront.

\section{Transitional Travelling Waves}\label{section7}
A \emph{transitional travelling wave }(or wavefront), in the present context, is defined as: \emph{a travelling wave of permanent form, propagating with constant propagation speed $v>0$.
The permanent wave form is nonnegative, and ahead of the wavefront the wave form approaches the equilibrium state $u_e=0$, whilst to the rear of the wavefront, the wave form either approaches the equilibrium state $u_e=1$ or approaches a non-trivial, positive periodic travelling wave (with the same propagation speed $v>0$).}
We anticipate that such structures will play a key role in the large-$t$ evolution of the solution to (IBVP), when $|x|=O(t)$ as $t\to\infty$.
This motivates the present consideration of transitional travelling waves. Following Section \ref{section6}, we can immediately eliminate the existence of transitional travelling waves of the second type described above, and we therefore concentrate on the existence of transitional travelling waves of the first type, which we will henceforth refer to as (TPTW) solutions to equation \eqref{1.3} with \eqref{1.9}. For a different treatment of permanent form transitional travelling waves to (\ref{eqn_NLFKPP0}), that lead up to the concept of generalised travelling waves, the reader is referred to \cite{Apreutesei2010} and \cite{VolpBio}  together with the references therein. In addition, in this context, numerical computations of transitional travelling waves have been considered in \cite{Nad}, relating particularly to such structures connecting two unstable steady states.

\begin{figure}
    \includegraphics[width=0.8\textwidth]{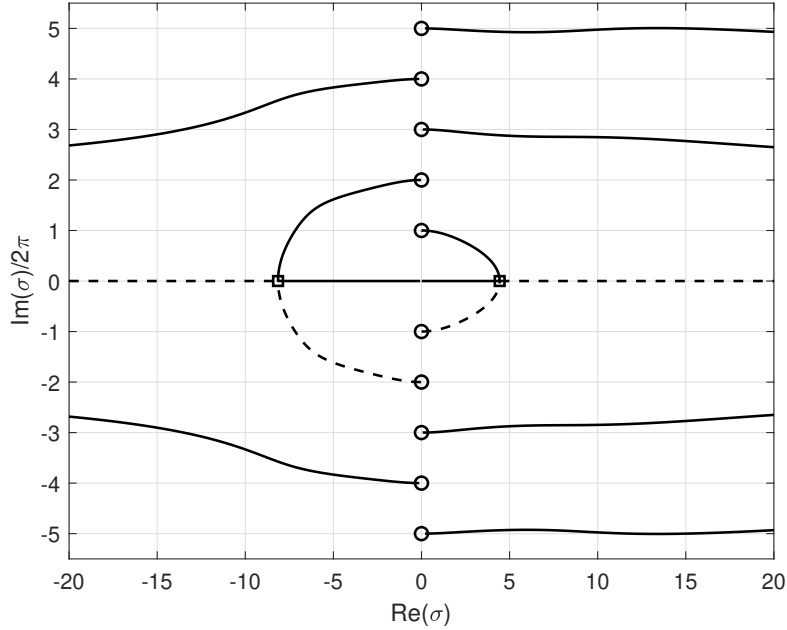}
    \caption{The locus of the eigenvalues $\sigma_n(D)$ in the complex plane. Note that the imaginary part has been scaled by a factor of $2 \pi$. These start from $\sigma_n = 2 \pi n i$ when $D=0$ (shown as circles) and move through the complex plane as $D$ increases. For $n = \pm 1$, $\pm 2$, $\sigma_n$ reaches the real axis at the points marked with a square, with $\sigma_1$, $\sigma_2 \to 0$ as $D \to \infty$, and $\sigma_{-1}\to\infty$, $\sigma_{-2} \to - \infty$ as $D \to \infty$. Both $\sigma_{-1}$ and $\sigma_{-2}$ are shown as broken lines.}
    \label{fig:7.1}
\end{figure}

We introduce the travelling coordinate $z=x-vt$ with $v>0$ being the constant propagation speed. We represent a (TPTW) solution with propagation speed $v>0$ as $u=u_T(z,v)$, with $u_T(\cdot,v)\in C^2(\R)$ satisfying,
\begin{align}
    \label{7.1}& D u_T''+v u_T'+u_T\left(1-\int_{z-\frac{1}{2}}^{z+\frac{1}{2}}u_T(s,v)ds\right)=0,\quad z\in\R
    \\\label{7.2}&u_T(z,v)\geq 0 \quad \forall~ z\in\R,
    \\\label{7.3}&u_{T}(z, v) \rightarrow\left\{\begin{array}{lll}
1 & \as & z \rightarrow-\infty \\
0 & \as & z \rightarrow+\infty.
\end{array}\right.
\end{align}
The problem \eqref{7.1}-\eqref{7.3} will be henceforth refered as (BVP).
We begin with some qualitative results.
\begin{itemize}
    \item[(T1)] Let $u_T(\cdot,v):\R\to\R$ be a (TPTW), then $u_T(\cdot,v)\in C^\omega(\R)$.
\end{itemize}
\begin{proof}
This follows by induction, following repeated differentiation through \eqref{7.1}, with the condition that $u_T(\cdot,v)\in C^2(\R)$.
\end{proof}
\begin{itemize}
    \item[(T2)] Let $u_T(\cdot,v):\R\to\R$ be a (TPTW), then $u_T(z,v)>0$ for all $z\in \R$. 
\end{itemize}
\begin{proof}
Suppose there is $z_0\in\R$ such that $u_T(z_0,v)=0$, then via \eqref{7.2} $u_T'(z_0,v)=0$. 
Thus, via induction on \eqref{7.1}, $u_T^{(n)}(z_0,v)=0$ for each $n=2,3,\dots$ . 
Since $u_T(\cdot,v)\in C^\omega(\R)$, via (T1), it then follows that $u_T(z,v)=0$ for all $z\in\R$, contradicting \eqref{7.3}. 
The result follows.
\end{proof}
\begin{itemize}
    \item[(T3)] The existence of a (TPTW) requires $v\geq 2\sqrt D$. 
\end{itemize}
\begin{proof}
Let $u_T(\cdot,v):\R\to\R$ be a (TPTW) with propagation speed $v>0$.
Via \eqref{7.1}-\eqref{7.3} and the linearization theorem (see for example \cite{Hale1977}), we must have,
\begin{align}
    \label{7.4}&D u_{T}^{\prime \prime}+v u_{T}^{\prime}+u_{T}=0, \quad z\gg1,
    \\\label{7.5}&u_T(z,v)\to0\,\as z\to\infty
    \\\label{7.6}&u_T(z,v)>0 ~\text{for} ~ z\gg1
\end{align}
with the strict inequality in \eqref{7.6} following from (T2).
As a consequence of \eqref{7.4} and \eqref{7.5} there exists constants $A_\infty$ and $B_\infty$, not both zero, such that,
\begin{align}
    \label{7.7}u_T(z,v)\sim A_\infty e^{\lambda_+(v)z}+B_\infty e^{\lambda_-(v)z},
\end{align}
with
\begin{align}
    \label{7.8}\lambda_{\pm}(v)=-\frac{1}{2 D}\left(v \pm\left(v^{2}-4 D\right)^{1 / 2}\right),
\end{align}
and the obvious modification when $v = 2\sqrt {D}.$ It then follows from \eqref{7.7} and \eqref{7.8} with \eqref{7.6}, that the propagation speed $v\geq 2 \sqrt{D}$, as required.
\end{proof}
In the remainder of this section we suppose that $u_T:\R\to\R$ is a (TPTW) which has the minimum possible propagation speed 
\begin{align}
    \label{7.9}v=2\sqrt{D}.
\end{align}
which we refer to as $u_T(z)$. Our intention is now to examine in detail the form of $u_T(z)$ as $z\to\ -\infty$, using the linearization theorem (see, for example, \cite{Hale1977}). We remark here that Fang and Zhao \cite{FZ} establish that there exists a critical value of $D$ such that $u_T(z)$ is monotone if and only if $D$ exceeds this critical value. We extend this here by obtaining details of this critical value, and establishing the existence of a second, and lower, critical value of $D$, for which the behaviour of $u_T(z)$ as $z\to -\infty$, although not globally monotone, is ultimately monotone exponential decay for $D$ above this second critical value, but becomes harmonically oscillating exponential decay for $D$ below this critical value.
First we write,
\begin{align}
    \label{7.10}u_T(z)=1+\ubar(z),
\end{align}
after which $\ubar (z)$ must satisfy
\begin{align}
    \label{7.11}&D\ubar''+2\sqrt{D}\ubar'-\int_{z-\frac{1}{2}}^{z+\frac{1}{2}}\ubar(s)ds=0,\quad(-z)\gg1, 
    \\\label{7.12}&\ubar(z)\to0\, \as z\to-\infty.
\end{align}
The form of $\ubar(z)$ as $z\to -\infty$ is determined by examining elementary solutions in the form
\begin{align}
    \label{7.13}\ubar (z)= e^{\sigma z}
\end{align}
with $\sigma\in\C$ to be determined, and $\op{Re}(\sigma)>0$, to satisfy condition \eqref{7.12} (both the real and imaginary parts of \eqref{7.13} provide real-valued solutions to the real linear equation \eqref{7.11}).
On substitution from \eqref{7.13} into \eqref{7.11} we find that $\sigma\in\C$  allows \eqref{7.13} to solve \eqref{7.11} and \eqref{7.12} if and only if 
\begin{align}
    \label{7.14}D\sigma ^2 + 2 \sqrt{D} \sigma -\frac{2}{\sigma} \sinh{\frac{1}{2}\sigma}=0, \text{ with } \op{Re}(\sigma)>0.
\end{align}
For each fixed $D> 0$, the transcendental equation \eqref{7.14} has a countably infinite number of roots in the complex plane, which we label as 
\begin{align}
    \label{7.15}\sigma_n (D)\in\C \text{ for each } n\in\Z\setminus\{0\}
\end{align}
and, in addition
\begin{align}
    \label{7.16}\sigma_n(\cdot)\in PC^1(\overline{\R}^+)\cap C(\overline{\R}^+).
\end{align}
with,
\begin{align}
    \label{7.17}\begin{split}
        &\op{Re}(\sigma_{\pm(2r-1)}(D))>0,\, r\in \N,
        \\&\op{Re}(\sigma_{\pm 2r}(D))<0,\, r\in \N.
    \end{split}
\end{align}
We next observe that,
\begin{align}
    \label{7.18} \sigma_n(D) = 2n\pi i -(-1)^n 8 n^2 \pi^2 \sqrt{D}+ O(D)\,\as D\to0^+
\end{align}
with $n\in\Z\setminus\{0\}$. Also,
\begin{align}
    \label{7.19}\begin{split}
        &\sigma_{\pm(2r-1)}(D)=2 \log \frac{1}{2}D \pm 4 (r-1)\pi i +o(1), 
        \\&\sigma_{\pm 2r}(D)=-2 \log \frac{1}{2}D \pm 4 (r-1)\pi i +o(1),
    \end{split}
\end{align}
as $D\to \infty$, with $r\in\N\setminus\{1\}$.
The full locus of $\sigma_n(D)$ is sketched in Figure~\ref{fig:7.1} for various values of $n$.
We observe from \eqref{7.14} that, in fact,
\begin{align}
    \label{7.20}\sigma_{-n}(D)=\bar{\sigma}_n(D)
\end{align}
for $n\in\Z\setminus \{\pm1,\pm2\}$.

We now consider the cases $n= \pm1,\pm2$.
We concentrate on $\sigma_{\pm1}(D)$.
For $0<D<D_+$, we have that $\sigma_{-1}(D)=\bar\sigma_{1}(D)$ and $\op{Im}(\sigma_1(D))>0$.
At $D=D_+$, $\sigma_{-1}(D_+)=\sigma_1(D_+)=\sigma_+\in\R^+$.
For $D>D_+$, $0<\sigma_1(D)<\sigma_+<\sigma_{-1}(D)$, with 
\begin{align}
    \label{7.21}\begin{split}
        &\sigma_1(D)\sim \frac{(\sqrt{2}-1)}{\sqrt{D}},
        \\&\sigma_{-1}(D)\sim 2\log\frac{1}{2}D ,
    \end{split}
\end{align}
as $D\to\infty$, as can be seen in Figure~\ref{fig:7.1}. In the above, $(\sigma_+,D_+)\in (\R^+)^2$ is the solution to the transcendental equations
\begin{align}
    \label{7.22}\begin{split}
        & D\sigma^3 +2 \sqrt{D}\sigma ^2 -2 \sinh \frac{1}{2}\sigma=0,
        \\& 3D \sigma ^2 + 4\sqrt{D} \sigma - \cosh\frac{1}{2} \sigma=0,
    \end{split}
\end{align}
which gives, via numerical approximation,
\begin{align}
    \label{7.23}(\sigma_+,D_+)\approx (4.437,2.824 \times 10^{-2}).
\end{align}

We are now able to interpret these results in relation to a (TPTW) with minimum propagation speed $v=2\sqrt D$.
In particular, \emph{translation invariance in $z$} can, \emph{and henceforth will be assumed to be}, fixed so that 
\begin{align}
    \label{7.24} u_T(z)\sim z e^{-\frac{1}{\sqrt{D}}z}\,\as z\to\infty
\end{align}
via \eqref{7.7} and \eqref{7.8}, whilst, for $0<D<D_+$, there are global constants $A_\infty\neq0$ and $\phi_\infty\in[0,2\pi)$ such that,
\begin{align}
    \label{7.25} u_T(z)\sim 1-A_\infty e^{a(D)z }\cos(b(D)z+\phi_\infty)\,\as z\to-\infty
\end{align}
with
\begin{align}
    \label{7.26} a(D)=\op{Re}(\sigma_1(D)), \quad  b(D)=\op{Im}(\sigma_1(D))
\end{align}
which are both positive.
In particular,
\begin{align}
    \label{7.27} & a(D) \sim\begin{cases}
8 \pi^{2} \sqrt{D} &\text{ as } D \rightarrow 0^{+}, \\
\sigma_{+}-O\left(\left(D_{+}-D\right)\right) &\text{ as } D \rightarrow D_{+}^{-}
\end{cases}
    \\\label{7.28} & b(D) \sim\begin{cases}
2 \pi+O(D) &\text{ as } D \rightarrow 0^{+}, \\
O((D_{+}-D)^{\frac{1}{2}}) &\text{ as } D \rightarrow D_+^-.
\end{cases}
\end{align}
However, for $D>D_+$, there is a global constant $A_\infty'\neq0$ such that \begin{align}
    \label{7.29} u_T(z)\sim 1-A_\infty' e^{\sigma_1(D)z}\,\as z\to-\infty
\end{align}
recalling that $0<\sigma_1(D)<\sigma_+$ and decreasing for $D\in (D_\infty,\infty)$, with
\begin{align}
    \label{7.30} \sigma_1(D)\sim \begin{cases}
\sigma_+-O((D-D_+)^{\frac{1}{2}})&\text{ as } D\to D^+_+  \\
\frac{(\sqrt2 -1)}{\sqrt D} & \text{ as } D\to\infty. 
\end{cases}
\end{align}

Principally we see that a (TPTW) with minimum speed is monotone as $z\to-\infty$ when $D\in(D_+,\infty)$, but has decaying oscillations as $z\to-\infty$ when $D\in(0,D_+)$. Numerical solutions of (IBVP) with $v=2\sqrt{D}$ and $D>0$ are shown in Figure~\ref{fig_TW}, and are consistent with the existence of a (unique by fixing translation invariance with $u_T(0)=\frac{1}{2}$) (TPTW) at minimum propagation speed $v = 2\sqrt{D}$, and that this is the travelling wave generated in the initial value problem (\ref{1.3}) to (\ref{1.5}) for $D>\Delta_1$. For $D<\Delta_1$, as discussed in Section~3, for moderate values of $D$, the minimum speed (TPTW) provides a good approximation to the solution of (IBVP) at the wavefront, which moves at speed $2\sqrt{D}$, ahead of a stationary periodic state. 

\section{Conclusions}\label{section8}
In this paper we have studied various aspects of the Cauchy problem for the nonlocal Fisher-KPP equation with a top hat kernel. We showed that the problem is globally well-posed and investigated the positive periodic steady state solutions that bifurcate from the uniform steady state $u=1$ at dimensionless diffusivity $D = \Delta_1 \approx 0.00297$. These have wavelength $\lambda$ and exist in a sequence of 'tongues' in $(\lambda, D)$ parameter space. As $D \to 0^+$, each of these periodic solutions has finite amplitude and wavelength, with $\lambda < 1$, and we constructed the asymptotic solution in detail for solutions in the tongue with largest wavelengths. 

We also investigated permanent form travelling wave solutions, which are known to exist for all wavespeeds $v$ with $v \geq 2 \sqrt{D}$, \cite{BNPR}. Numerical solutions suggest that, from localised initial conditions with finite support, a pair of diverging travelling wavefronts with minimum wavespeed, $2 \sqrt{D}$, is generated, and that for $D \geq \Delta_1$, the solution asymptotes to the stable uniform state, $u=1$, between the wavefronts. For $0<D<\Delta_1$ and moderately small ($\sim 10^{-3}$), a periodic static steady state, with wavelength around $0.7$, is generated between the wavefronts. The wavelength selection mechanism is not clear, but it is notable that $0.7$ is approximately the linearly most unstable wavelength of the uniform steady state $u=1$, and that this instability is not dispersive. However, as $D$ becomes very small ($\sim 10^{-6}$), this wavelength drifts towards $0.5$, the minimum available wavelength in the first tongue, with spike formation initiating in the region ahead of the wavefront where $u$ is exponentially small, and a rational argument for this mechanism has been presented in section 2.

In the next paper in this series, \cite{BN2023_1} we will report on the global bifurcation structure for the problem set on a finite domain with various boundary conditions. Oscillatory solutions are also generated, with a structure similar to those discussed in the present paper, but with nontrivial dependence on the form of the boundary conditions.

An interesting extension is to the Cauchy problem in higher dimensions, for example with the two dimensional kernel
\[
\phi(r) = \left\{
\begin{array}{cc}
\frac{4}{\pi} & \mbox{for $r<\frac{1}{2}$,}\\
0 & \mbox{otherwise,}
\end{array}\right.
\]
in polar coordinates $(r,\theta)$. The numerical solution of initial value problems in higher spatial dimensions is significantly more challenging than in one dimension. We would expect that spatially-localised solutions would be generated for small enough $D$, and that they would be amenable to asymptotic analysis.

\section*{Acknowledgements}

N.M. Ladas was supported by the EPSRC Doctoral Training Partnership, Reference: EP/R513167/1.

The authors would like to thank the referees for very helpful comments, which have led to a significant improvement in the paper.

\section*{Competing Interests} The authors declare none.

\printbibliography

@article{VolpBio,
author = {Banerjee, M. and Kuznetsov, M. and Udovenko, O. and Volpert, V.},
journal = {Acta Biotheoretica},
title = {Nonlocal Reaction-Diffusion Equations in Biomedical Applications},
year = {2022},
volume = {70}
}

@article{Nad,
author = {Nadin, G. and Perthame, B. and Tang, M.},
journal = {C. R. Acad. Sci. Paris, Ser. 1},
title = {Can a travelling wave connect two unstable states, The case of the nonlocal Fisher equation},
year = {2011},
volume = {349}
}

@article{Fay,
author = {Faye, G. and Holzer, M. },
journal = {J. Diff. Eqns.},
title = {Modulated traveling fronts for a nonlocal Fisher-KPP equation : A dynamical systems approach},
year = {2015},
volume = {258}
}

@article{Hamel2014,
author = {Hamel, F. and Ryzhik, L.},
doi = {10.1088/0951-7715/27/11/2735},
issn = {13616544},
journal = {Nonlinearity},
title = {{On the nonlocal Fisher-KPP equation: Steady states, spreading speed and global bounds}},
year = {2014}
}

@book{VanDykeMilton1975Pmif,
publisher = {Parabolic Press},
booktitle = {Perturbation methods in fluid mechanics},
isbn = {0915760010},
year = {1975},
title = {Perturbation methods in fluid mechanics / by Milton Van Dyke.},
edition = {Annotated ed.},
language = {eng},
address = {Stanford, Calif.},
author = {Van Dyke, M.},
lccn = {76-351560},
}

@article{Apreutesei2010,
author = {Apreutesei, N. and Bessonov, N. and Volpert, V. and Vougalter, V.},
doi = {10.3934/dcdsb.2010.13.537},
issn = {15313492},
journal = {Discrete and Continuous Dynamical Systems - Series B},
number = {3},
pages = {537--557},
title = {{Spatial structures and generalized travelling waves for an integro-differential equation}},
volume = {13},
year = {2010}
}

@article{Volpert2009,
author = {Volpert, V. and Petrovskii, S.},
doi = {10.1016/j.plrev.2009.10.002},
issn = {15710645},
journal = {Physics of Life Reviews},
number = {4},
pages = {267--310},
publisher = {Elsevier B.V.},
title = {{Reaction-diffusion waves in biology}},
url = {http://dx.doi.org/10.1016/j.plrev.2009.10.002},
volume = {6},
year = {2009}
}

@book{Kavallarisbook,
address = {Cham},
author = {Kavallaris, N. I. and Suzuki, T.},
doi = {10.1007/978-3-319-67944-0},
isbn = {978-3-319-67942-6},
pages = {300},
publisher = {Springer International Publishing},
series = {Mathematics for Industry},
title = {{Non-Local Partial Differential Equations for Engineering and Biology}},
url = {http://link.springer.com/10.1007/978-3-319-67944-0},
volume = {31},
year = {2018}
}

@book{Hale1977,
author = {Hale, J. K.},
booktitle = {In Applied Mathematical Sciences},
edition = {3},
issn = {00765392},
title = {{Introduction to Functional Differential Equations}},
volume = {99},
year = {1977}
}

@article{BNPR,
Author  = "H. Berestycki and G. Nadin and B. Perthame and L. Ryzhik",
Title   = "The non-local {F}isher-{KPP} equation: travelling waves and steady states",
Journal = "Nonlinearity",
Volume  = "22",
Pages   = "2813-2844",
Year    = "2009"
}

@article{JBNL,
	doi = {10.1088/1361-6544/ab6f4f},
	url = {https://doi.org/10.1088/1361-6544/ab6f4f},
	year = 2020,
	publisher = {{IOP} Publishing},
	volume = {33},
	number = {5},
	pages = {2106--2142},
	author = {J. Billingham},
	title = {Slow travelling wave solutions of the nonlocal {F}isher-{KPP} equation},
	journal = {Nonlinearity}
}

@article{JB1,
  title={Dynamics of a strongly nonlocal reaction-diffusion population model},
  author={J. Billingham},
  journal={Nonlinearity},
  year={2004},
  volume={17},
  pages={313-346}
}

@misc{BN2023_1,
      title={The evolution problem for the 1D nonlocal Fisher-KPP equation with a top hat kernel. Part 2. The Cauchy problem on a finite interval}, 
      author={D. J. Needham and J. Billingham},
      year={2023},
      eprint={2304.10935},
      archivePrefix={arXiv},
      primaryClass={math.AP}
}

@book{Coombes,
Author    = " S. Coombes and P. Graben and R. Potthast and J. Wright",
Title     = "Neural Fields: Theory and Applications",
Publisher = "Springer",
Year      = "2014"
}

@Article{Furter1989,
author="Furter, J.
and Grinfeld, M.",
title="Local vs. non-local interactions in population dynamics",
journal="Journal of Mathematical Biology",
year="1989",
volume="27",
number="1",
pages="65--80"
}

@Article{Gourley2000,
author="Gourley, S.A.",
title="Travelling front solutions of a nonlocal Fisher equation",
journal="Journal of Mathematical Biology",
year="2000",
volume="41",
number="3",
pages="272--284"
}

@Article{Van,
author={van Saarloos, W.},
title="Front propagation into unstable states",
journal="Physics Reports",
year="2003",
volume="386",
number="2",
pages="29--222"
}

@Article{BHR,
author={E. Bouin and C. Henderson and L. Ryzhik},
title="The Bramson delay in non local Fisher-KPP equation",
journal="Ann. I. H. Poincare",
year="2020",
volume="37",
pages="51--77"
}

@article{coville_dupaigne_2007, title={On a non-local equation arising in population dynamics}, volume={137}, DOI={10.1017/S0308210504000721}, number={4}, journal={Proceedings of the Royal Society of Edinburgh: Section A Mathematics}, publisher={Royal Society of Edinburgh Scotland Foundation}, author={Coville, J. and Dupaigne, L.}, year={2007}, pages={727–755}}

@article{perthame2007, 
title="Concentration in the Nonlocal Fisher Equation: the Hamilton-Jacobi Limit", 
number="4", 
journal="Mathematical Modelling of Natural Phenomena",
author="Perthame, B. and Génieys, S.", 
year="2007", 
pages="135–151"}

@incollection{VV,
  author      = "V. Volpert and V. Vougalter",
  title       = "Emergence and Propagation of Patterns in Nonlocal Reaction-Diffusion Equations Arising in the Theory of Speciation",
  editor      = "M. Lewis and Ph. Maini and S. Petrovsky",
  booktitle   = "Dispersal, individual movement and spatial ecology",
  publisher   = "Springer",
  year        = 2013,
  pages       = "331-353"}

@article{Britton,
author = "Britton, {N. F.}",
year = "1990",
month = "12",
volume = "50",
pages = "1663--1688",
journal = "SIAM Journal on Applied Mathematics",
number = "6",
}

@Article{Gierer1972,
author="Gierer, A.
and Meinhardt, H.",
title="A theory of biological pattern formation",
journal="Kybernetik",
year="1972",
volume="12",
number="1",
pages="30--39"
}

@Inbook{Volpert2014,
author="Volpert, V.",
title="Nonlocal Reaction-diffusion Equations",
bookTitle="Elliptic Partial Differential Equations: Volume 2: Reaction-Diffusion Equations",
year="2014",
publisher="Springer Basel",
address="Basel",
pages="521--626"
}

@article{FZ,
    author = {Fang, J. and Zhao, X-Q.},
    title = "{Monotone wavefronts of the nonlocal Fisher-KPP equation}",
    journal = {Nonliearity},
    volume = {24},
    pages = {3043-3054},
    year = {2011},
}

@article{Surulescu,
author = {Li, J. and Chen, L. and Surulescu, C.},
title = {Global boundedness, hair trigger effect, and pattern formation driven by the parametrization of a nonlocal {F}isher-{KPP} problem},
journal = {Journal of Differential Equations},
volume = {269},
number = {11},
pages = {9090-9122},
year = {2020},
issn = {0022-0396},
doi = {https://doi.org/10.1016/j.jde.2020.06.039},
url = {https://www.sciencedirect.com/science/article/pii/S0022039620303600}
}

@book{LNB,
author = {J. A Leach and D. J. Needham},
title={Matched {A}symptotic {E}xpansions in {R}eaction-{D}iffusion {T}heory},
year={2003},
publisher={Springer Monographs in Mathematics},
address={London},
}

@book{KUZ,
author = {Y. A. Kuznetsov},
title={Elements of {A}pplied {B}ifurcation {T}heory},
year={1995},
publisher={Springer},
address={New York},
}

@article{GVA,
    author = {S. Genieys and V. Volpert and P. Auger},
    title = "{Patterns and waves for a model in population dynamics with nonlocal consumption of resources}",
    journal = {Mathematical Modelling of Natural Phenomena},
    volume = {1},
    pages = {65-82},
    year = {2006},
}

\end{document}